\documentclass[letter,12pt]{article}

\usepackage{preamble}
\usepackage[english]{babel}
\usepackage[utf8]{inputenc}
\usepackage{tikz-cd}
\usetikzlibrary{shapes,arrows}
\tikzstyle{block} = [rectangle, draw,
    text width=7em, text centered, rounded corners,draw=black, minimum height=2.5em,text centered]
\tikzstyle{line} = [draw, -latex', thick]
\tikzset{node distance = 2cm}
\theoremstyle{plain}
\newtheorem*{rep@theorem}{\rep@title}
\newcommand{\newreptheorem}[2]{%
\newenvironment{rep#1}[1]{%
 \def\rep@title{#2 \ref{##1}}%
 \begin{rep@theorem}}%
 {\end{rep@theorem}}}
\makeatother
\newtheorem{thmx}{Theorem}

\newreptheorem{theorem}{Theorem}
\author{{\scshape David Crnčević}~~and~~{\scshape Felipe Hernández}~~and~~{\scshape Kevin Rizk}\\~~and~~{\scshape Khunpob Sereesuchart}~~and~~{\scshape Ran Tao}\\}
\date{\small \today}
\title{\bfseries On the multiplicative independence between $n$ and $\lfloor \alpha n\rfloor$}

\begin{document}

\maketitle
\begin{abstract}
In this article we investigate different forms of multiplicative independence between the sequences $n$ and $\lfloor n\alpha\rfloor$ for irrational $\alpha$. Our main theorem shows that for a large class of arithmetic functions $a,b\colon\N\to\C$ the sequences $(a(n))_{n\in \N}$ and $(b(\lfloor \alpha n\rfloor))_{n\in \N}$ are asymptotically uncorrelated. This new theorem is then applied to prove a $2$-dimensional version of the \Erdos-Kac theorem, asserting that the sequences $(\omega(n))_{n\in \N}$ and $(\omega(\lfloor\alpha n\rfloor))_{n\in \N}$ behave as independent normally distributed random variables with mean $\log\log n$ and standard deviation $\sqrt{\log\log n}$. Our main result also implies a variation on Chowla's Conjecture asserting that the logarithmic average of $(\lambda(n) \lambda ( \lfloor \alpha n\rfloor))_{n\in \N}$ tends to $0$. 

\end{abstract}
\tableofcontents
\thispagestyle{empty}

\section{Introduction}
For a set $A\subseteq \N$ we say $A$ has \emph{density} $d(A)$ if $d(A):=\lim_{n\to \infty}\frac{|A\cap\{1,\ldots,N\}|}{N}$ exists. Even though $d$ is not a probability measure, since it is not countably additive, it can be used to measure the size of infinite sets of integers. For instance, if we define $A_p$ to be the set of positive integers divisible by $p\in \primes$, then $A_p$ has density $d(A_p)=1/p$.\\

It is believed that the sequences $(n)_{n\in \N}$ and $(\lfloor \alpha n\rfloor)_{n\in \N}$, for $\alpha\in \R\setminus\Q$, have independent multiplicative structure. The following classical result of G. L. Watson \cite{Watson53} from 1953 supports this claim, showing that the proportion of natural numbers $n\in \N$ such that $n$ and $\lfloor \alpha n\rfloor$ are coprime is $6/\pi^2$, which is the same probability of two arbitrary integers being relatively prime.

\begin{Theorem}[see {\cite{Watson53}}]\label{Watson}
For every $\alpha\in \R\setminus\Q$:
$$d(\{ n\in \N \mid (n,\lfloor \alpha n\rfloor )=1\})= \frac{6}{\pi^2}, $$
were $(a,b)$ denotes the greatest common divisor between $a,b\in \N$.
\end{Theorem}

In 1917, Hardy and Ramanujan \cite{HR17} showed the following theorem which describes the asymptotic behaviour of the function $\omega(n)$ which counts distinct prime factors of $n$.

 \begin{Theorem}[see {\cite{HR17}}]\label{Hardy-Ramanujan}
Let $\epsilon>0$, then for almost all $n\in \N$
$$|\omega(n)-\log\log n|<{(\log\log n)}^{\frac12 +\varepsilon}.$$
\end{Theorem}

In $1940$ P. Erd\"{o}s and M. Kac  proved their well-known result about the Gaussian law associated to the function $\omega$, stated here.
\begin{Theorem}[see {\cite{EK40}}]\label{ErdosKac}
For all $a,b \in \R$ where $a < b$,
$$d\left(\left\{n\in\N : a \leq \frac{\omega(n)-\log\log{n}}{(\log\log{n})^{1/2}}\leq b\right\}\right)=\frac{1}{\sqrt{2\pi}}\int_{a}^b e^{-t^2/2 }dt.$$
\end{Theorem}

Later in 2007, William D. Banks and Igor E. Shparlinski showed that one can replace $n$ in the \Erdos{}-Kac theorem by any Beatty sequence (see \cref{beattysequence}).

\begin{Theorem}[see {\cite[Theorem 4]{BeattySeq}}]\label{ErdosKacBS}
Let $\alpha\in \R$ be a positive irrational real number, then
$$d\left(\left\{n\in\N : a \leq \frac{\omega(\lfloor \alpha n\rfloor)-\log\log{n}}{(\log\log{n})^{1/2}}\leq b\right\}\right)=\frac{1}{\sqrt{2\pi}}\int_{a}^b e^{-t^2/2 }dt.$$
\end{Theorem}

Our first result unifies the theorem of Banks and Shparlinski with the original result of \Erdos{}-Kac, following the heuristic that the sequences $n$ and $\lfloor \alpha n\rfloor$, for $\alpha$ irrational, behave independently.

\begin{thmx}\label{ErdosKacGeneral}
Let $\alpha\in \R$ be a positive irrational real number and $a,b,c,d$ real numbers. Then
\begin{align*}
d\left(\left\{n\in \N : a\leq \frac{\omega(\lfloor \alpha n\rfloor)-\log\log{n}}{(\log\log{n})^{1/2}}\leq b\right\}\cap \left\{n\in \N : c\leq \frac{\omega(n)-\log\log{n}}{(\log\log{n})^{1/2}}\leq d\right\}\right)\\
=\Bigl(\frac{1}{\sqrt{2\pi}}\int_{a}^be^{-t^2/2 }dt\Bigr)\Bigl(\frac{1}{\sqrt{2\pi}}\int_{c}^d e^{-t^2/2 }dt\Bigr).    
\end{align*}
\end{thmx}

For our second result, we took inspiration from a long-standing conjecture in number theory attributed to Chowla. 
\begin{Conjecture}[Chowla's conjecture [see {\cite[Problem 57]{Chowla66} ,\cite{Tao16}]}]\label{Chowlas-conjecture}
Let $k\geq 1$ be a natural number.
\begin{enumerate}
    \item (Chowla's conjecture) If $h_1,\cdots,h_k\in \Z$ are distinct integers, then
    $$\lim_{N\to\infty} \E_{n\leq N} \lambda(n+h_1)\cdots\lambda(n+h_k)=0. $$
    \item (Logarithmically averaged Chowla's conjecture) If $h_1,\cdots,h_k\in \Z$ are distinct integers, then
    $$\lim_{N\to\infty} \E_{n\leq N}^{\text{log}} \lambda(n+h_1)\cdots\lambda(n+h_k)=0, $$ 
\end{enumerate} 
where $\lambda$ is the Liouville function (see \cref{def_lio}) and $\E_{n \leq N}$, $\E_{n \leq N}^{\log}$ are standard and logarithmic averages (see \cref{defaverages}).
\end{Conjecture}
\begin{Remark}
The logarithmically averaged Chowla's conjecture was first stated by Tao in \cite{Tao16}. Furthermore, Chowla's conjecture (with \Cesaro{} averages) implies the logarithmically averaged version.
\end{Remark}

Chowla's conjecture is known to be true only for $k=1$, which case is elementally equivalent to the prime number theorem. In addition, due to recent developments, we know that the logarithmically averaged Chowla's conjecture holds for $k= 2$ \cite{Tao16} and for odd values of $k$ \cite{Tao2018}. Inspired by this recent breakthrough, we established the following result.
\begin{thmx}\label{LiouvilleTheorem}
If $\alpha$ is an irrational positive real number, then 
$$\lim_{N\rightarrow \infty} \E_{n\leq N}^{\log} \lambda(n)\lambda(\lfloor \alpha n\rfloor)=0,$$
where $[N]:=\{1,\cdots,N\}$. 
\end{thmx}

Both Theorem \ref{ErdosKacGeneral} and Theorem \ref{LiouvilleTheorem} are consequences of a more general theorem about a class of bounded sequences $a:\N\to \C$ which we call \emph{\textbf{BMAI}} (\emph{Bounded Multiplicative Almost Invariant}). Among other things, this class of sequences includes all multiplicative functions which don't take the value $0$; for more details and further examples, see \cref{BMAI}. 

\begin{thmx}\label{Main_Theorem}
Let $a: \N \to \C$ and $b: \N \to \mathbb{S}^1$ be BMAI sequences such that:
\begin{enumerate}[topsep=4pt, parsep=4pt]
    \item $\lim_{H \to \infty} \limsup_{N \to \infty} \E_{n\leq N} | \E_{h\leq H} a(n+h)| = 0$; and
    \item for all $\beta\in\R\setminus \Q$ and $h\in \N$ we have
    $$\lim_{N \to \infty}  \E_{n\leq N} a(n) \overline{a(n+h)} e(n \beta) =0 ~~~ ( \text{resp. }\E^{\log}_{n\leq N} \bullet).$$
\end{enumerate}
Then, for any irrational $\alpha \in \R$,
\begin{equation}
    \lim_{N \to \infty} \E_{n\leq N} b(n) a(\lfloor \alpha n \rfloor) = 0. ~~~( \text{resp. replacing } \E_{n\leq N} \text{ with }\E^{\log}_{n\leq N} \bullet)
\end{equation} 
\end{thmx}

To derive \cref{ErdosKacGeneral} from \cref{Main_Theorem}, we take $a(n)=G(\frac{\omega(n)-\log\log{n}}{(\log\log{n})^{1/2}})$ and $b(n)=F(\frac{\omega(n)-\log\log{n}}{(\log\log{n})^{1/2}})$ for certain functions $F,G:\N\to \C$, and employ a generalization of the \Erdos{}-Kac theorem on short intervals to verify that the hypothesis of \cref{Main_Theorem} are satisfied for this choice of $a$ and $b$. The details of this derivation are given in \cref{sec3} and \cref{sec4}. \\

Likewise, for the proof of \cref{LiouvilleTheorem} we take $a(n)=b(n)=\lambda(n)$ in \cref{Main_Theorem} and use \cref{Matomaki-Radziwill} and \cref{LiouvilleConditions}. For the details of the proof of \cref{LiouvilleTheorem}, we refer the reader to \cref{sec3} and \cref{sec4} as well.\\

Finally, since there are many intermediate steps involved in the proof of \cref{Main_Theorem}, we include here a diagram illustrating the path to our main results.
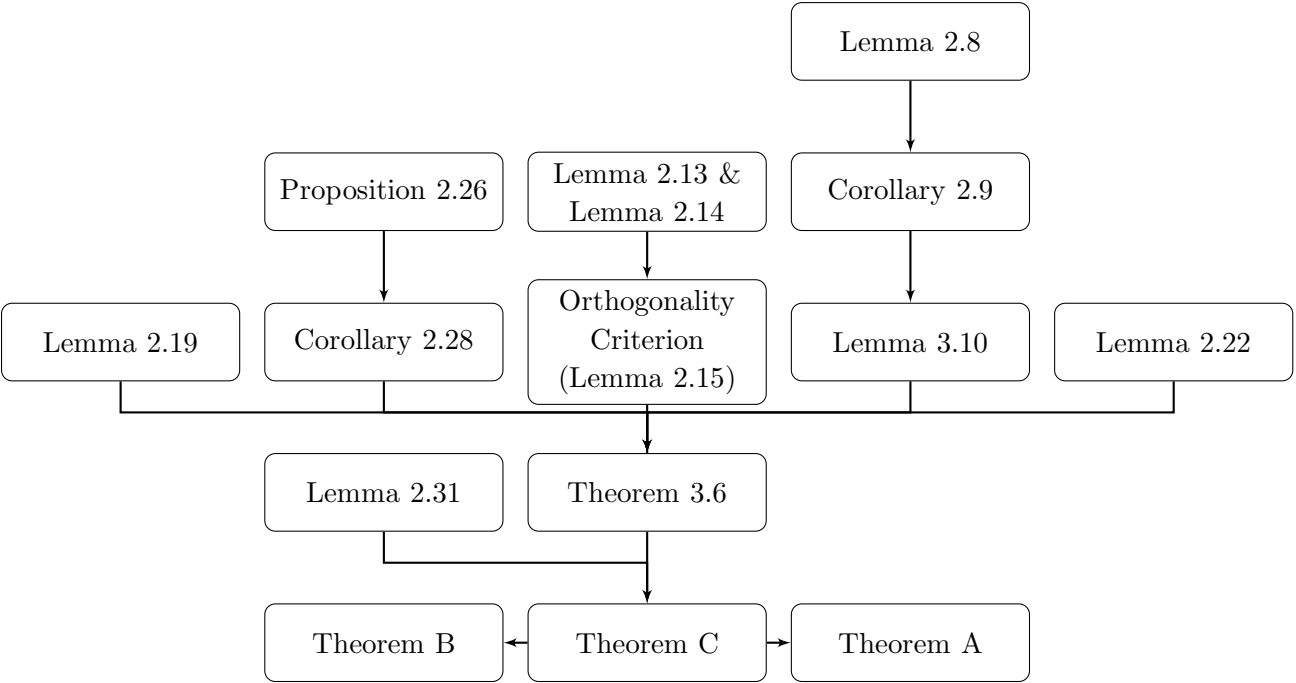
\begin{center}
\begin{tikzpicture}
 \node [block] (LemmasforOC) {\small \cref{Turan_dual} \& \cref{lemma_log_inx_chang} };

   \node [block, below of= LemmasforOC] (OCriterion) {\small Orthogonality Criterion (\cref{orthogonality})};
   \node [block, left of=OCriterion, node distance=3.5cm] (Corollary228) {\small \cref{Lemma_BMAI}};
   \node [block, left of=Corollary228, node distance=3.5cm] (Lemma219) {\small \cref{APApproximation}};
   \node [block, right of=OCriterion, node distance=3.5cm] (Lemma310) {\small \cref{previousLemma}};
    \node [block, right of=Lemma310, node distance=3.5cm] (Lemma222) {\small \cref{DpiAP}};

 \node [block, above of= Lemma310] (Corollary29) {\small \cref{Lemma_zero_Ces_Log}};
 \node [block, above of=Corollary29] (Lemma28) {\small \cref{Limsup-log-average}};

 \node [block, above of= Corollary228] (Proposition226) {\small \cref{BMAItwoprimes}};

 \node [block, below of=OCriterion]  (Theorem36)  {\small \cref{technicalTheorem}};
\node  [block, left of= Theorem36, node distance=3.5cm] (Lemma231) {\small \cref{1->2}}; 
  \node [block, below of=Theorem36]  (TheoremC)  {\small \cref{Main_Theorem}};

  \node  [block, right of= TheoremC, node distance=3.5cm] (TheoremA) {\small \cref{ErdosKacGeneral}};
 \node  [block, left of= TheoremC, node distance=3.5cm] (TheoremB) {\small \cref{LiouvilleTheorem}};

\draw [-, thick]  (OCriterion.south) --  (Theorem36);

\draw[|-,-|,-latex', thick,] (Lemma310.south) |-+(0,-1em)-| (Theorem36.north); 

\draw[|-,-|,-latex', thick,] (Lemma222.south) |-+(0,-1em)-| (Theorem36.north);   

\draw[|-,-|,-latex', thick,] (Lemma219.south) |-+(0,-1em)-| (Theorem36.north);   

\draw[|-,-|,-latex', thick,] (Corollary228.south) |-+(0,-1em)-| (Theorem36.north);

\draw [line]  (Lemma28) -- (Corollary29);
\draw [line]  (Corollary29) -- (Lemma310);
\draw [line]  (LemmasforOC) -- (OCriterion);
\draw [line]  (Proposition226) -- (Corollary228);

\draw[-latex', thick,] (Lemma231.south) |-+(0,-1em)-| (TheoremC.north);  
\draw [line]  (Theorem36) --  (TheoremC);
\draw [line]  (TheoremC) --  (TheoremA);
\draw [line]  (TheoremC) -- (TheoremB);
 
\end{tikzpicture}

\end{center}

\subsection*{Acknowledgements}
We would like to thank the Bernoulli Center for funding the Young Researchers in Mathematics Summer Program, which allowed us to develop such results. We would also like to thank Florian Richter and Joel Moreira for their guidance throughout the  program and for making this project possible. Additionally, we owe a special thank you to Noah Kravitz and Katharine Woo for providing detailed comments and suggestions. 
\section{Preliminaries}
\subsection{Notation}
First, we introduce some classical functions in number theory.

\begin{Definition}[Prime omega functions]
We define the prime omega functions $\Omega$,  $\omega:\N\to \N$ as follows.
For $n\in \N$, we write $n=p_1^{a_1}\cdots p_k^{a_k}$ where $p_1,\dots,p_k$ are distinct primes numbers and $a_j\in\N$ for $j = 1, \dots, k$. We define
$$\omega(n)\coloneqq k.$$ 
In other words, $\omega(n)$ is the number of distinct prime divisors of $n$. We also define
$$\Omega(n)\coloneqq \sum_{j=1}^k a_j,$$
which is the numbers of prime divisors with multiplicity.
\end{Definition}
\begin{Definition}[Liouville function]
\label{def_lio}
The classical Liouville function is defined as
$$\lambda(n) \coloneqq (-1)^{\Omega(n)}.$$
\end{Definition}


Throughout the paper we will use $e(x) \coloneqq e^{2\pi i x}$, $x\in \R$, for the complex exponential function. For a real number $\alpha$, we will denote $\lfloor\alpha\rfloor$ the greatest integer less than or equal to $\alpha$, and $\{\alpha\}=\alpha-\lfloor\alpha\rfloor$ the fractional part. We now define the notion of a Beatty sequence.

\begin{Definition}[Beatty Sequence]\label{beattysequence}
For $\alpha\in \R_+$ the corresponding Beatty sequence is the sequence of integers defined by
$$\mathcal{B}_{\alpha}=(\lfloor \alpha n\rfloor)_{n\in \N}. $$
\end{Definition}

Finally, we introduce some classical definitions to measure how fast a sequence goes to $0$ and how bounded the growth rate of a sequence is.
\begin{Definition}
For a function $g:\N\to \R$, which is strictly positive for all large enough values of $n$, we use the classic notation $f(N)=o_{N\to\infty}(g(N))$ to denote a function $f:\N\to \R$ such that 
$$\lim_{N\to\infty} \frac{f(N)}{g(N)}=0. $$
Additionally, we denote $f(N)=O{\bigl (}g(N){\bigr )}$, if there exists a positive real number $C$ and an integer $M\in \N$ such that
$$|f(N)|\leq Cg(N), \forall N\geq M. $$
\end{Definition}

\subsection{Averaging}
In our proofs, we will sometimes use expectation notation to denote averages. In particular, we define here the \Cesaro~and the logarithmic averages.
\begin{Definition}\label{defaverages} 
    For a function $f : \N \to \C$ and  $N,S\in \N$, we define its
    \begin{itemize}
        \item \emph{average} as 
        $$\E_{n \leq N}(f) \coloneqq \frac{1}{N}\sum_{n=1}^N f(n),$$
        \item \textit{logarithmic average} as 
        $$\E_{n \leq N}^{\log}(f) \coloneqq \frac{1}{\log N}\sum_{n=1}^N \frac{f(n)}{n},$$
        \item \textit{logarithmic expectation over the primes} as 
    \end{itemize}
    \begin{equation*}
        \E_{p \leq S}^{\log} (g) \coloneqq \frac{1}{\log \log S}\sum_{p \leq S} \frac{1}{p} g(p),
    \end{equation*}
    where $p$ only takes prime values less than or equal to $S$.
\end{Definition}
\begin{Remark}
When the limit $\lim_{N\to\infty} \E_{n\leq N}(f)$ exists, it is called the \emph{\Cesaro{} average} of the sequence $f:\N\to \C$. In this case, the  limit $\lim_{N\to\infty} \E_{n\leq N}^{log}(f)$ exists as well.
\end{Remark}
We will start by stating and proving results which relates the 
\Cesaro~average and the logarithmic average.\\

\begin{Lemma} \label{log_ces}
We have for $f: \N \to \C$ bounded, 
\begin{align*}
    \E_{n \leq N}^{\log} f(n) = \E_{M \leq N}^{\log} \E_{n \leq M} f(n) + o_{N \to \infty}(1).
\end{align*}
\end{Lemma}

\begin{proof}
By expanding the sum $\E_{M \leq N}^{\log} \E_{n \leq M} f(n)$, we get 
\begin{align*}
    \E_{M \leq N}^{\log} \E_{n \leq M} f(n) &= \frac{1}{\log N} \sum_{M = 1}^ N \frac{1}{M^2} \sum_{n = 1} ^ M f(n)  \\
    &= \frac{1}{\log N}  \sum_{n = 1}^N f(n) \sum_{M = n}^N \frac{1}{M^2}.
\end{align*}
Using 
\begin{align*}
\sum_{M = n}^N \frac{1}{M^2} = \frac{1}{n} + \frac{1}{N} + O(1/n^2),
\end{align*}
and that $f$ is bounded, we get:
\begin{align*}
\E_{M \leq N}^{\log} \E_{n \leq M} f(n) &= 
\frac{1}{\log N} \sum_{n=1}^N \Big( \frac{f(n)}{n}  + \frac{f(n)}{N} + O(1/n^2)  \Big)\\
&= \E_{n \leq N}^{\log} f(n) + \frac{1}{N \log N} \sum_{n = 1}^N O(1) + \frac{1}{\log N} \sum_{n = 1}^N O(1/n^2) \\
&= \E_{n \leq N}^{\log} f(n) + o_{N \to \infty}(1).
\end{align*}
\end{proof}

\begin{Lemma}\label{Limsup-log-average}
Let $f:\N\to \C$ a bounded function. Then
$$\limsup_{N\to\infty}\E^{\log}_{n \leq N} |f(n)|\leq\limsup_{N\to\infty}\E_{n \leq N} |f(n)|.$$
\end{Lemma}

\begin{proof}
For each $N\in \N$, define 
$$M_N=\text{argmax}\{\E_{n\leq M} |f(n)| ~\big| \log\log{N} \leq M\leq N\}.$$
 By \cref{log_ces} we have that:
   \begin{align*}
     \E^{\log}_{n \leq N} |f(n)| &=  \E_{M \leq N}^{\log} \E_{n \leq M} |f(n)| +o_{N\to \infty}(1)\\
     &=   \frac{1}{\log{N}}\sum_{M<\log\log{N}} \frac{1}{M} \E_{n \leq M} |f(n)|+ \frac{1}{\log{N}}\sum_{M=\log\log{N}}^{N}\frac{1}{M} \E_{n \leq M} |f(n)| +o_{N\to \infty}(1)\\
     &\leq \frac{\Vert f\Vert _{\ell^\infty} \log\log{N}}{\log{N}}+\Big( \frac{1}{\log{N}}\sum_{M=\log\log{N}}^{N}\frac{1}{M}\Big)\Big( \E_{n \leq M_N} |f(n)|\Big) +o_{N\to \infty}(1),
 \end{align*}
applying $\limsup_N$ yields 
\begin{align*}
      \limsup_{N \to \infty } \E^{\log}_{n \leq N} |f(n)|&\leq \limsup_{N} \E_{n \leq M_N} |f(n)|\leq \limsup_{N}\E_{n \leq N} |f(n)|,
\end{align*}
as desired.
\end{proof}
The following corollary will be useful later on.
\begin{Corollary}\label{Lemma_zero_Ces_Log}
    Let $S \subseteq \N$ such that $|S| = \infty$ and  let $f: \N \times S \to \C$ such that for all $t \in S$, $f(\cdot,t)$ is a bounded function. We have that
    $$\limsup_{T \to \infty, T \in S}\limsup_{N \to \infty } \E^{\log}_{n \leq N} |f(n,T) | \leq  \limsup_{T \to \infty, T \in S}\limsup_{N \to \infty } \E_{n \leq N} |f(n,T) |.$$
\end{Corollary}

 The rest of this subsection will contain some useful known results.
 
\begin{Lemma}[see {\cite[Theorem 2.4]{BerMor16}}]\label{corput}
    Let $(u_n)_{n\in \N}$ be a bounded sequence in $\mathbbm{C}$. If for all $h\in \N$ 
    \begin{equation}
        \frac{1}{N}\sum_{n=1}^N u_{n+h} \overline{u_n} \rightarrow 0,
    \end{equation}
    then 
    \begin{equation}
        \frac{1}{N}\sum_{n=1}^N u_n \rightarrow 0.
    \end{equation}
\end{Lemma}

The following lemmas, scattered throughout the literature, will allow us to transition between \Cesaro{} averages and logarithmic expectations over the primes.

\begin{Lemma}\label{Turan-Kub_ineq}
For $S, N \in \N$ with $N > S$, let $L(S) = \sum_{p \leq S} \frac{1}{p}$, we have 
\begin{align} \label{eq_tk_ces}
    \E_{n \leq N} \big| \sum_{p \leq S, p |n } 1 - L(S)\big|^2 \leq E(S,N), 
\end{align}
with $E(S,N) = O(\log \log S) + \frac{1}{N}O(S^2)$. We also have 
\begin{align} \label{eq_tk_log}
    \E_{n\leq N}^{\log} \big| \sum_{p \leq S, p |n } 1 - L(S)\big|^2 \leq E'(S,N),
\end{align}
with $E'(S,N) = O(\log \log S) + \frac{1}{\log(N)}O(S^2)$.
\end{Lemma}

\begin{proof} \cref{eq_tk_ces}  follows from \cite[Lemma 4.1]{Elliott79} applied to the additive arithmetic function $f(n) := \sum_{p \leq S, p|n} 1$. \\\\
Now let's prove \cref{eq_tk_log}.\\
First, we estimate $\sum_{n \leq N} \sum_{p \leq S, p | n } \frac{1}{n}$: 
\begin{align*}
\sum_{p \leq S} \sum_{n \leq N, p | n} \frac{1}{n}  &=  \sum_{p \leq S} \sum_{n \leq N/p} \frac{1}{np} = \sum_{p \leq S} \frac{1}{p} \sum_{n \leq N/p} \frac{1}{n}\\ 
&= \sum_{p \leq S} \frac{1}{p}(\log(N) - \log(p) + O(1)) = L(S)\log(N) + O(S). 
\end{align*}
Second, we estimate $ \sum_{n\leq N} \sum_{p,q \leq S, p,q | n } \frac{1}{n} $:
\begin{align*}
\sum_{n\leq N} \sum_{p,q \leq S, p,q | n } \frac{1}{n} & =  \sum_{p \leq S} \sum_{q \leq S} \sum_{n \leq N, p | n, q|n } \frac{1}{n} \\
&=  \sum_{p \leq S} \sum_{n \leq N, p | n} \frac{1}{n} + \sum_{p \leq S} \sum_{q \leq S, q \neq p} (\sum_{n \leq N/pq} \frac{1}{pqn} )  \\
& \leq L(S)\log(N) + O(S) + L(S)^2 \log(N) + O(S^2).
\end{align*}
Finally, we estimate $\frac{1}{\log N }\sum_{n \leq N} \frac{1}{n} \big| \sum_{p \leq S, p |n } 1 - L(S)\big|^2$:
\begin{align*}
    &\frac{1}{\log N}\sum_{n \leq N} \frac{\big| \sum_{p \leq S, p |n } 1 - L(S)\big|^2}{n} \\
    &\tab= \frac{1}{\log N}\sum_{n\leq N} \sum_{p,q \leq S, p,q | n } \frac{1}{n} -\Bigl( 2L(S) \frac{1}{\log N}\sum_{n \leq N}\sum_{p \leq S, p |n } \frac{1}{n} \Bigr)+ L(S)^2.
\end{align*}
Now using our estimations, we get
\begin{align*}
    \frac{1}{\log N }\sum_{n \leq N} \frac{\big| \sum_{p \leq S, p |n } 1 - L(S)\big|^2}{n} &\leq  L(S) + \frac{1}{\log N }O(S)  +\frac{1}{\log N}O(S^2)\\
    &\tab+ L(S)^2 -2L(S)^2 +\frac{L(S)}{\log N}O(S) + L(S)^2\\
    &= L(S)  +\frac{1}{\log N }O(S^2),
\end{align*}
and as $L(S) = O(\log \log S)$, we get the desired result.\\
\end{proof}
\begin{Remark}
    We note that the above lemma still holds when dropping the condition $N > S$ by observing that when $p>N$, $\sum_{n \leq \frac{N}{p}} \frac{1}{n} = 0$.
\end{Remark}
\begin{Lemma}\label{Turan_dual}
For any bounded function $f : \N \to \C$, we have that
    \begin{equation}
        \frac{1}{N}\sum_{n=1}^N f(n) = \mathbb{E}_{p \leq S}^{\log} \frac{1}{N/p} \sum_{n = 1}^{N/p} f(pn) + \epsilon(S,N),
    \end{equation}
with $\limsup_{S \to \infty}\limsup_{N \to \infty}\epsilon(S,N) = 0$.
We also have that
    \begin{equation}\label{lemma2.9log}
        \E_{n \leq N}^{\log} f(n) = \mathbb{E}_{p \leq S}^{\log} \E_{n\leq N/p}^{\log}f(pn) + \epsilon'(S,N),
    \end{equation}
with $\limsup_{S \to \infty}\limsup_{N \to \infty}\epsilon'(S,N) = 0$.
\end{Lemma}

\begin{proof}
Let $L(S) = \sum_{p \leq S} \frac{1}{p}$,   \cref{Turan-Kub_ineq} gives
\begin{align*}
    \frac{1}{N}\sum_{n \leq N} \big| \sum_{p \leq S, p |n } 1 - L(S)\big|^2 \leq E(S,N).
\end{align*}
By Cauchy-Schwarz, and as $f$ is bounded we get
\begin{align*}
    \big | \frac{1}{N}\sum_{n \leq N} f(n)\big( \sum_{p \leq S, p |n } 1 - L(S)\big) \big |  \leq C E(S,N)^{\frac{1}{2}},
\end{align*}
for some constant $C$ independent of $N$ and $S$. Now dividing by $L(S)$, we get
\begin{align*}
    \big | \frac{1}{N}\sum_{n \leq N} f(n) - \frac{1}{L(S)} \frac{1}{N} \sum_{n \leq N} \sum_{p \leq S, p |n } f(n) \big | \leq  \frac{1}{L(S)} C E(S,N)^{\frac{1}{2}},
\end{align*}
where the second term on the left can be written as
\begin{align*}
    \frac{1}{L(S)N} \sum_{n \leq N} \sum_{p \leq S, p |n } f(n)  &= \frac{1}{L(S)} \sum_{p \leq S} \frac{1}{N} \sum_{n \leq N , p | n} f(n) = \frac{1}{L(S)}\sum_{p \leq S} \frac{1}{p} \frac{1}{N/p} \sum_{n \leq N/p } f(pn) \\
     &= \mathbb{E}_{p \leq S}^{\log} \frac{1}{N/p} \sum_{n = 1}^{N/p} f(pn).
\end{align*}
The only thing left to show is that the error term goes to zero as desired. Indeed,
\begin{align*}
    \limsup_{S \to \infty} \limsup_{N \to \infty} \frac{1}{L(S)} C E(S,N)^{\frac{1}{2}} &= \limsup_{S \to \infty} \frac{1}{L(S)}\left (\limsup_{N \to \infty} O( \log \log S) + \frac{1}{N} O\left(S^2\right) \right )^{\frac{1}{2}} \\
    &= \limsup_{S \to \infty } O\left((\log \log S)^{\frac{-1}{2}}\right) = 0 ,
\end{align*}
where we used that $L(S) = O( \log \log S)$ to conclude.\\

The proof for \eqref{lemma2.9log} is similar.

\end{proof}

\begin{Lemma} \label{lemma_log_inx_chang} 
For any bounded $f : \N \to \C $ and for any $k,N \in \N$, we have
\begin{align*}
    \E_{n \leq N}^{\log} f(n) = \E_{n \leq N/k}^{\log} f(n) + E(N,k)
\end{align*}
such that $\limsup_{N \to \infty}E(N,k) = 0$.
\end{Lemma}
\begin{proof}
By using the definition of logarithmic average, we have
\begin{align*}
    E(N,k) = \bigg(\frac{1}{\log N} - \frac{1}{\log (N/k)} \bigg )  \sum_{n = 1}^{N/k} \frac{f(n)}{n} + \frac{1}{\log N}  \sum_{n = N/k}^{N} \frac{f(n)}{n}. 
\end{align*} 
We will start by showing  that the first part of the sum tends to zero as $N$ goes to infinity. Indeed, we have
\begin{equation}\label{zero-expression}
    \bigg(\frac{1}{\log N} - \frac{1}{\log (N/k)} \bigg )  \sum_{n = 1}^{N/k} \frac{f(n)}{n} = \frac{1}{\log N} \sum_{n = 1}^{N/k} \frac{f(n)}{n} \big( \frac{\log{k}}{\log N - \log k} \big ) ,
\end{equation}
but as $f$ is bounded, we have that $\frac{1}{\log N} \sum_{n = 1}^{N/k} \frac{f(n)}{n}$ is bounded, and as $\big( \frac{\log{k}}{\log N - \log k} \big )$ tends to $0$ as $N$ tends to infinity, we get that \eqref{zero-expression} goes to zero as $N \to \infty$.\\

As for the second  part of the sum, as $f$ is bounded, we have that $|\limsup_{N \to \infty} \sum_{n = N/k}^{N} \frac{f(n)}{n} |  \leq \max_{n \in \N} | f(n) | \log k $ and so we indeed get that
\begin{align*}
    \limsup_{N \to \infty} \frac{1}{\log N}  \sum_{n = N/k}^{N} \frac{f(n)}{n} = 0.
\end{align*}
\end{proof}

Using \cref{Turan_dual}, we can derive the following orthogonality criterion, which has appeared in \cite{Katai86}, \cite{DD82} and \cite{BSZ13}. 

\begin{Lemma}\label{orthogonality}
There exists $C$ such that for every bounded function $f:\N \to \C$ and $N,S\in \N$
\begin{equation}\label{Florian-inequality}
    | \E_{n \leq N} f(n) | \leq C \sum_{l = 1}^{\log_2 S} \frac{1}{l \log \log S} \left(\max_{p,q \in [2^l,2^{l+1})} \left|\E_{n \leq \frac{N}{\max \{p,q\}}} f(pn) \overline{f(qn)}\right|\right)^{1/2} + \epsilon(S,N),
\end{equation}
where $\limsup_{ S \to \infty } \limsup_{N \to \infty } \epsilon(S,N) = 0$.

We also have the following logarithmic version of the inequality:
\begin{equation}\label{Florian-inequality_log}
    \left| \E_{n \leq N}^{\log} f(n) \right| \leq \left( \E_{p \leq S}^{\log}\E_{q \leq S}^{\log}  \left | \E^{\log}_{n \leq N} f(pn)\overline{f(qn)} \right | \right)^{1/2}+ \epsilon'(S,N)
\end{equation}
where $\limsup_{ S \to \infty } \limsup_{N \to \infty } \epsilon'(S,N) = 0$.

\end{Lemma}
\begin{proof}
Let $A_l:=[2^l,2^{l+1})$. 
\begin{align*}
    \left| \E_{n \leq N} f(n) \right| &= \left|\E_{p \leq S}^{\log} \E_{n \leq N/p} f(pn) \right| + \epsilon(S,N)\text{ by \cref{Turan_dual}} \\
    &\leq \sum_{l =1}^{{\log_2 s} } \left|\E_{p \leq S}^{\log} \one_{A_l}(p) \left(\frac{2^l}{p} \E_{n \leq N / 2^l} \one_{pn \leq N} f(pn)\right)\right| + \epsilon(S,N) \\
    &= \sum_{l =1}^{\log_2 s } \left|\E_{n \leq N / 2^l} \E_{p \leq S}^{\log} \one_{A_l}(p) \left(\frac{2^l}{p} \one_{pn \leq N} f(pn)\right)\right| + \epsilon(S,N) \\
    &\leq \sum_{l =1}^{\log_2 s } \left(\E_{n \leq N / 2^l} \left|\E_{p \leq S}^{\log} \one_{A_l}(p) \left(\frac{2^l}{p} \one_{pn \leq N} f(pn)\right)\right|^2\right)^{1/2} + \epsilon(S,N) \\
    &= \sum_{l =1}^{\log_2 s } \left(\E_{n \leq N / 2^l} \E_{p,q \leq S}^{\log} \one_{A_l}(p)\one_{A_l}(q) \left(\frac{4^l}{pq} \one_{pn \leq N} \one_{qn \leq N} f(pn)\overline{f(qn)}\right)\right)^{1/2} + \epsilon(S,N)  \\ 
    &\leq 2 \sum_{l =1}^{\log_2 s } \left( \E_{p,q \leq S}^{\log} \one_{A_l}(p)\one_{A_l}(q) \left( \max_{p,q \in A_l}| \E_{n \leq N / 2^l} \one_{\leq N}(pn) \one_{\leq N}(qn) f(pn)\overline{f(qn)}|\right)\right)^{1/2} \\
    &\tab+ \epsilon(S,N).
\end{align*}
 We can expand 
 $$\E_{p,q \leq S}^{\log} \one_{A_l}(p)\one_{A_l}(q) = \frac{1}{(\log \log S)^2} \sum_{p \in A_l} \frac{1}{p} \sum_{q \in A_l} \frac{1}{q}. $$ 
 Using the prime number theorem, we also have 
 $$\sum_{p \in A_l} \frac{1}{p} = \sum_{q \in A_l} \frac{1}{q} \leq \sum_{q \in A_l} \frac{1}{2^l} \leq \theta \frac{1}{l},$$ 
 giving \eqref{Florian-inequality}. \\
 Now let's prove \eqref{Florian-inequality_log}: 
 \begin{align*}
    \left| \E_{n \leq N}^{\log} f(n) \right| &= \left|\E_{p \leq S}^{\log} \E^{\log}_{n \leq N/p} f(pn) \right| + \epsilon(S,N)\text{ by \cref{Turan_dual}} \\
    &\leq \left|\E^{\log}_{n \leq N}\E_{p \leq S}^{\log}  f(pn) \right| +  \left|\E_{p \leq S}^{\log} E(N,p)  \right| + \epsilon(S,N) \text{ by \cref{lemma_log_inx_chang}} \\
    &\leq \big( \E^{\log}_{n \leq N}\left|\E_{p \leq S}^{\log}  f(pn) \right|^2\big)^{1/2} +  \left|\E_{p \leq S}^{\log} E(N,p)  \right| + \epsilon(S,N) \\
    &= \big( \E^{\log}_{n \leq N}\E_{p \leq S}^{\log}\E_{q \leq S}^{\log}  f(pn)\overline{f(qn)}\big)^{1/2}+ \epsilon'(S,N)\\ 
    &\leq \big( \E_{p \leq S}^{\log}\E_{q \leq S}^{\log}  \left | \E^{\log}_{n \leq N} f(pn)\overline{f(qn)} \right |)^{1/2}+ \epsilon'(S,N).
\end{align*}

\end{proof}

Next is a powerful theorem from \cite{Tanaka}, formalizing the independent behavior of polynomials with regard to the function $\Omega$.
\begin{Theorem}[see {\cite[Theorem 1]{Tanaka}}]\label{tanaka}
    Let $\{f_i: i\leq k\}$ be a finite family of pairwise relatively prime non-constant polynomials, where $f_i$ is a product of $r_i \geq 1$ irreducible polynomials. \\
    Let 
    \begin{equation}
        u_i(n) \coloneqq \frac{\Omega(f_i(n))-r_i \log\log n}{\sqrt{r_i \log \log n}}, 
    \end{equation} 
    then for $E \subset \R^k$ Jordan measurable, 
    \begin{equation}
        \lim_{N \rightarrow \infty} \frac{1}{N} \# \{n \leq N: (u_i(n))_i \in E \} = \frac{1}{2\pi^\frac{k}{2}} \int_E \exp(\frac{-1}{2} \sum_{i=1}^k x_i^2 ) dx_1 \dots dx_k.
    \end{equation}
\end{Theorem}
\begin{Remark}
\cref{tanaka} is also true when $\Omega$ is replaced by $\omega$, proven in the same paper. 
\end{Remark}

\subsection{Almost Periodicity}
We now recall the notion of almost periodic function, which plays an important role in the sequel. Given $g:\N\to \C$, we define its $1$-norm as
$$\|g\|_1 \coloneqq \limsup_{N \to \infty} \E_{n\leq N}|g(n)|. $$
\begin{Definition}
Let $f:\N\to \C$ be an arithmetic function. We say that $f$ is \emph{almost periodic} (sometimes also referred to as \emph{Besicovitch} \emph{almost periodic} in the literature) if for each $\epsilon>0$ there is some linear combination $h$ over $\C$ of exponential functions $e(\alpha n)$, $\alpha\in \R$, such that $\|f-h\|_1<\epsilon$.
\end{Definition}

\begin{Lemma}\label{APApproximation}
Let $f:\N\rightarrow \C$ and $g:\R\rightarrow \C$ be two functions such that $f=g$ on $\N$ and $g$ is a Riemann integrable bounded periodic function with period $\alpha\in \R\setminus \Q$. Then, for all $\epsilon>0$, $\exists M\in \N$ and $c_k\in\C$
$$P(x)=\sum_{k=0}^M c_k  e(kx/\alpha),$$
such that $\Vert f-P\Vert_1\leq \epsilon$. In particular, $f$ is almost periodic.
\end{Lemma}
\begin{Remark}
If $\alpha\in \Q$, the result still holds. In fact, in that case $f$ is periodic.
\end{Remark}
\begin{proof}
Consider the function $h:\R\to \C$ defined as $h(x)=g(x\alpha)$. Note that $h$ is a $1$-periodic function and piece-wise continuous. As $h$ is bounded, it can be seen as a function of $L^2([0,1])$. Thus, there exists 
$c_k:=\hat{h}(k)=\frac{1}{2\pi}\int_{-\pi}^\pi h(x) e^{ikx} dx,$

where $\hat h(k)$ is the $k$th Fourier coefficient of $h$, such that $\Vert h-\sum_{k=0}^Mc_ke(kx)\Vert _{L^1[0,1]}\to_M 0$.\\ 

Take $M$ big enough such that $\Vert h-\sum_{k=0}^Mc_ke(kx)\Vert _1\leq \epsilon$, and define $P(x)=\sum_{k=0}^Mc_ke(kx)$. Since $f(n)=h(n/\alpha)$, we have that 
$$f(n)- \sum_{k=0}^Mc_ke(kn/\alpha)= h(n/\alpha)-P(n/\alpha),$$
thus
\begin{equation}\label{eq2.11}
    \frac{1}{N}\sum_{n=1}^N|f(n)- \sum_{k=0}^Mc_ke(kn/\alpha)|=\frac{1}{N}\sum_{n=1}^N|h-P|(n\alpha^{-1}).
\end{equation}
As $(n\alpha^{-1})_n$ is equidistributed in $\T$ (as $\alpha^{-1}$ is irrational), and $|h-P|$ is a piece-wise continuous function on the torus, by the Weyl equidistribution criteria we have that 
$$\frac{1}{N}\sum_{n=1}^N|h-P|(n\alpha^{-1})\xrightarrow[N\to \infty]{} \int_{\T} |h-P| d\mu =\Vert h-P\Vert _1\leq \epsilon. $$
Therefore, applying limsup in \eqref{eq2.11}, we conclude that $\Vert f-P\Vert _1\leq \epsilon$. 

\end{proof}

\begin{Lemma}\label{IndIsAP}
If $\alpha\in \R$, then $\one_{\lfloor \alpha \Z\rfloor}:\N\to \C $ is an almost periodic function, which can be approximated with linear combination of exponential functions $e(\alpha n)$ with $\alpha\in \R\setminus\Q$. 
\end{Lemma}
\begin{proof}
First, we note that for $m\in \Z$,
\begin{equation}
\begin{split}
m \in \lfloor \alpha \Z\rfloor &\iff \exists  n\in \Z \text{ such that } m \leq \alpha n < m+1 \\
&\iff \exists  n\in \Z \text{ such that } \alpha n -1 < m \leq \alpha n.
\end{split}
\end{equation}
Consider the set 
$$R_\alpha \coloneqq \bigcup_{n\in \Z} (\alpha n-1,\alpha n].$$
We note that $\one_{\lfloor \alpha \Z\rfloor}=\one_{R_{\alpha}}$ on $\N$. Treating $\one_{R_{\alpha}}$ as a function over $\R$, we have that it is a Riemann integrable bounded periodic function with period $\alpha\in \R$. Therefore, \cref{APApproximation} yields the conclusion. 
\end{proof}
\begin{Lemma}\label{DpiAP}
For $p\in \primes$, fix $0\leq i<p$, $\alpha\in \R\setminus \Q$ and let 
$$D_{p,i} \coloneqq \{n\in \lfloor \alpha \N\rfloor: n=\lfloor m\alpha \rfloor \text{ and } \lfloor p \alpha m\rfloor = p\lfloor \alpha m \rfloor +i\}. $$
There is a Riemann integrable $\alpha$-periodic function $f$ over the reals such that $f=\one_{D_{p,i}}$ when restricted to $\N$.

\end{Lemma}
\begin{proof}
First, for any $m \in \N$,
\begin{equation}
    \begin{split}
        \lfloor p\alpha m\rfloor = p\lfloor \alpha m \rfloor + i &\iff \{\alpha m\} \in \left[\frac{i}{p}, \frac{i+1}{p}\right)\\
        &\iff \alpha m - \frac{i+1}{p} < \lfloor \alpha m\rfloor \leq \alpha m - \frac{i}{p}.
    \end{split}
\end{equation}
Therefore, we define the set
$$D_{p,i}':=\bigcup_{m\in \Z} \left(\alpha m - \frac{i+1}{p}, \alpha m -\frac{i}{p}\right],$$
and note that $\one_{D_{p,i}}=\one_{D_{p,i}'}$ on $\N$. Furthermore, $f:=\one_{D_{p,i}'}$ is a Riemann integrable function, which is also $\alpha$-periodic.
\end{proof}

\subsection{BMAI sequences}
The goal of this subsection is to introduce the notion of BMAI sequences and list some properties. We let $\primes$ denote the set of prime numbers. 
\begin{Definition}\label{BMAI}
A bounded sequence $a : \N \to \C$ is called \emph{BMAI} if there is a bounded sequence $(\theta_p)_{p\in \primes}\subseteq \C\setminus\{0\}$ such that
\begin{equation}
    \limsup_{N \to \infty } \E_{n \leq N}| a(pn) - \theta_p a (n)|  = 0.
\end{equation}
\end{Definition}
The idea behind this definition is that the limiting behaviour of a BMAI sequence approximates that of a multiplicative function.

\begin{Remark}
If $a : \N \to \C$ is a BMAI sequence we also have by \cref{Lemma_zero_Ces_Log} using $S = \P$ and for $T\in S$, $f(n,T) = a(Tn) - \theta_T a(n)$, that
\begin{equation*}
    \limsup_{N \to \infty } \E_{n \leq N}^{\log}| a(pn) - \theta_p a (n)|  = 0.
\end{equation*}
\end{Remark}
\begin{Example}
The following are BMAI functions.
\begin{itemize}
    \item The Liouville function $\lambda:\N\to \C$ with $\theta_p=-1$, $\forall p\in \primes$.
    \item The function $\psi_F(n) \coloneqq F\left(\frac{\Omega(n) - \log\log n}{\sqrt{\log\log n}}\right)$, for a bounded continuous $F: \R \to \C$ with $\theta_p=1$, $\forall p\in \primes$. 
    \item Any bounded multiplicative arithmetic function $f$, with $\theta_p=f(p)$ nonzero, $\forall p\in \primes$. 
\end{itemize}
\end{Example}
We prove some basic properties. 

\begin{Proposition}\label{BMAItwoprimes}
Let $a,b:\N\to \C$ BMAI sequences. Then for primes $p,q\in \primes$ we have 
$$\limsup_N \E_{n \leq N} \left(\left|a(pn)b(qn)-\theta_p^a\theta_q^b a(n)b(n)\right|\right) = 0,$$
and also
$$\limsup_N \E_{n \leq N}^{\log} \left(\left|a(pn)b(qn)-\theta_p^a\theta_q^b a(n)b(n)\right|\right) =0.$$
\end{Proposition}
\begin{proof}
We have the following inequality
\begin{align*}
    \left|a(pn)b(qn)- \theta_p^a\theta_q^ba(n)b(n)\right| &\leq \left|a(pn)b(qn)-\theta_p^a a(n) b(qn)\right|+\left|\theta_p^a a(n) b(qn)- \theta_p^a\theta_q^ba(n)b(n)\right|\\
    &\leq \Vert b\Vert_{\ell^\infty}\left|a(pn)-\theta_p^aa(n)\right|+ \Vert a\Vert_{\ell^\infty}\Vert(\theta_p^a)_p\Vert_{\ell^\infty} \left|b(pn) -\theta_q^bb(n)\right|
\end{align*}
which gives
$$\limsup_N\left| a(pn)b(pn)- \theta_p^a\theta_p^ba(n)b(n)\right|=0.$$
\end{proof}
\begin{Corollary}
BMAI sequences are closed under coordinate-wise multiplication, i.e. if $(a(n))_n$ and $(b(n))_n$ are BMAI then so is $(a(n)b(n))_n$.
\end{Corollary}
\begin{proof}
Taking $p=q$, \cref{BMAItwoprimes} gives 
\begin{align*}
    \limsup_{N \to \infty } \E_{n \leq N}| (a\cdot b)(pn) - \theta_p^a\theta_p^b (a\cdot b)(n)|  =0.
\end{align*}
\end{proof}

\begin{Corollary} \label{Lemma_BMAI}
For a BMAI sequence $a : \N \to \C$, 
$$\limsup_N\E_{n \leq N} |a(pn)\overline{a(qn)} -  \theta_p \overline{\theta_q}|a(n)|^2|=0,$$ and
$$\limsup_N\E_{n \leq N}^{\log} |a(pn)\overline{a(qn)} -  \theta_p \overline{\theta_q}|a(n)|^2|=0,$$
for all $p, q \in\primes$. 
\end{Corollary}  
\begin{proof}
Take $b=\overline{a}$ in Proposition \ref{BMAItwoprimes}.
\end{proof}
\begin{Proposition}\label{BMAI_over_progressions}
Let $a:\N\to \C$ be a BMAI sequence and $f:\N\to \C$ a bounded function. Then for any $d_1,d_2\in \N$
$$\limsup_{N\to \infty}\E_{n\leq N}|a\Big(p(d_1n+d_2)\Big)f(n) -\theta_p a(d_1n+d_2)f(n)| =0,  $$
and 
$$\limsup_{N\to \infty}\E_{n\leq N}^{\log}|a\Big(p(d_1n+d_2)\Big)f(n) -\theta_p a(d_1n+d_2)f(n)| =0.  $$
\end{Proposition}
\begin{proof}
Without loss of generality, we can take $f\equiv 1$ since 
$$\E_{n\leq N}|a(p(d_1n+d_2))f(n) -\theta_p a(d_1n+d_2)f(n)|\leq \Vert f\Vert _{\ell^\infty} \E_{n\leq N}|a(p(d_1n+d_2))-\theta_p a(d_1n+d_2)|. $$
Note that 
\begin{align*}
    \E_{n\leq N}|a(p(d_1n+d_2))-\theta_p a(d_1n+d_2)|& =   \frac{1}{N} \sum_{m = 1}^{d_1N+d_2} |a(pm)-\theta_p a(m)| \one_{ d_1\N+d_2}(m)\\
    &=  \left(d_1+\frac{d_2}{N}\right) \E_{m\leq d_1N+d_2}|a(pm)-\theta_p a(m)| \one_{(d_1\N+d_2)}(m)\\
    &\leq \left(d_1+\frac{d_2}{N}\right) \E_{m\leq d_1N+d_2}|a(pm)-\theta_p a(m)|.
\end{align*}
Taking $\limsup$, we derive 
$$\limsup_{N\to \infty}\E_{n\leq N}|a(p(d_1n+d_2))f(n) -\theta_p a(d_1n+d_2)f(n)| =0. $$
\end{proof}
\begin{Proposition}
Let $a: \N\to \C$ be a BMAI sequence, then there exists a (actual) multiplicative function $\theta:\N\to \C$ such that for every $m\in \N$
$$\limsup_{N\to \infty}\E_{n\leq N} |a(mn)-\theta(m)a(n)| =0. $$
Additionally, $||a||_1\neq 0$ if and only if $\theta$ is unique.
\end{Proposition}
\begin{proof}
Let $(\theta_p)_{p\in \P}$ be a sequence given by $a$ being BMAI. We define $\theta:\N\to \C$, as the multiplicative function defined by $\theta(p)=\theta_p$, $\forall p\in \P$. \\

We prove the statement using induction over $\Omega(m)$. The definition of a BMAI sequence yields the case $\Omega(m)=1$. We assume that the statement holds for all $m\in \N$ with $\Omega(m) \leq k$. Let $m\in \N$ with $\Omega(m)=k$ and $p\in \P$. We use \cref{BMAI_over_progressions} with $d_1=m$, $d_2=0$ and $f\equiv 1$ to obtain
\begin{equation}\label{eq2.15}
 \limsup_{N\to \infty}\E_{n\leq N} | a( pmn)-\theta_p a(mn)|=0.   
\end{equation}
Using \eqref{eq2.15} and the inductive hypothesis: 
\begin{align*}
    &\limsup_{N\to \infty}\E_{n\leq N} |a(pmn)-\theta(pm)a(n)|\\
    &\tab\leq  \limsup_{N\to \infty}\E_{n\leq N} |a(pmn)-\theta_pa(mn)|+ |\theta_p|\E_{n\leq N} |a(mn)-\theta(m)a(n)|\\
    &\tab=0,
\end{align*}
and as every $n\in \N$ with $\Omega(n)=k+1$ can be written as $n=pm$ for $p\in \P$ with $p\mid n$ and $m=n/p$, thus concluding existence.\\

For uniqueness, take $p\in \P$ and let $\theta_p'$ such that
$$\limsup_{N\to \infty} |a(pn)-\theta_p'a(n)|=0. $$
Note that
$$\limsup_{N\to\infty} \E| \theta_p a(n)-\theta_p'a(n)|\leq \limsup_{N\to \infty }  |a(pn)-\theta_p'a(n)|+ |a(pn)-\theta_pa(n)| =0. $$ 
 Therefore 
 $$|\theta_p-\theta_{p}'|\cdot ||a||_1=0, \forall p\in \P, $$
 so if $||a||_1\neq 0$ then $\theta_p=\theta_p'$, $\forall p\in \P$. For the other direction, if $||a||_1=0$, then $\theta_p\in \C$ can be any value. In fact, notice that
 \begin{align*}
      \limsup_{N\to \infty} \E_{n\leq N} |a(pn)-\theta_p a(n)| &\leq \limsup_{N\to \infty}\E_{n\leq N} |a(pn)| +|\theta_p|\limsup_{N\to \infty}\E_{n\leq N}|a(n)|\\
      &\leq p||a||_1 + |\theta_p|\cdot ||a||_1=0,
 \end{align*}
which doesn't depend on $\theta_p$.
\end{proof}

\begin{Lemma}\label{1->2} 
Let $a:\N\to\C$ be a BMAI sequence such that for all natural numbers $h\in \N$ and irrational $\beta$, 
$$\lim_{N\to\infty} \E_{n \leq N} a(n)\overline{a(n+h)}e(\beta n)=0~~~ ( \text{resp. }\E^{\log}_{n\leq N} \bullet).$$
Then, $\forall p,q\in \primes$, $\beta\in \R\setminus\Q$ and $h_1,h_2\in \N$ with $qh_1-ph_2\neq 0$,
$$\limsup_{N\to\infty} \E_{n \leq N}  a(pn+h_1)\overline{a(qn+h_2)}e(\beta n)=0.$$
\end{Lemma}

\begin{proof}
We write out the proof for \Cesaro{} averages only, as the logarithmic case is analogous. Without loss of generality, assume $ph_2-qh_1>0$. Then
\begin{align*}
\lim_{N\to\infty} \E_{n\leq N} a(n+qh_1)\overline{a(n+ph_2)}e(\beta n) e(\beta qh_1) =\lim_{N\to\infty} \E_{n\leq N} a(n)\overline{a(n+ph_2-qh_1)}e(\beta n) =0.
\end{align*}
Thus
\begin{equation}\label{2.14}
    \lim_{N\to\infty} \E_{n\leq N}a(n+qh_1)\overline{a(n+ph_2)}e(\beta n) =0,
\end{equation}
for every $\beta\in \R\setminus\Q$ and $p,q\in \primes$ and $h_1,h_2\in \N$ with $qh_1-ph_2\neq 0$. On the other hand, using \cref{BMAI_over_progressions} twice, we get
$$\limsup_{N\to\infty}\E_{n\leq N}\Big|\overline{\theta}_p\theta_q a(pn+h_1)\overline{a(qn+h_2)}e(\beta n)- a(pqn+qh_1)\overline{a(pqn+ph_2)}e(\beta n)\Big| =0.$$

Hence
\begin{align*}
&\limsup_{N\rightarrow \infty} \Big|\E_{n\leq N} \overline{\theta}_p\theta_qa(pn+h_1)\overline{a(qn+h_2)}e(\beta n)\Big|\\
&\tab=\limsup_{N\rightarrow \infty} \Big|\E_{n\leq N}a(pqn+qh_1)\overline{a(pqn+ph_2)}e(\beta n)\Big| \\
&\tab=\limsup_{N\rightarrow \infty} \Big|\frac{1}{N} \sum_{n=1}^{Npq}a(n+qh_1)\overline{a(n+ph_2)}e(\frac{\beta}{pq} n)\one_{pq\N}(n)\Big|.
\end{align*}
Since $\one_{pq\N}(n)=\frac{1}{pq}\sum_{k=1}^{pq} e(nk/pq)$, we have 

\begin{align*}
 &\limsup_{N\rightarrow \infty} \Big|\E_{n \leq N}  a(pn+h_1)\overline{a(qn+h_2)}e(\beta n)\Big|\\
 &\tab= \Big|\frac{1}{\overline{\theta}_p\theta_q}\Big|\cdot\Bigl(\sum_{k=1}^{pq} \lim_{N\rightarrow \infty} \frac{1}{N} \Big|\sum_{n=1}^{Npq} a(n+qh_1)\overline{a(n+ph_2)}e(\frac{\beta+k}{pq}n)\Big|\Bigr)\\
 &\tab=0,
\end{align*}
 in which we used equation \eqref{2.14} in the last equality to conclude.
\end{proof}

\section{Main Result}\label{sec3}
\subsection{Basic Case}
The following result is well known, as discussed in \cite{jiang_lu_wang_2020}, for multiplicative sequences. We generalize it to all BMAI sequences.
\begin{Lemma} \label{Lemma_BMAI_and_exp}
Let $a : \N \to \C$ be a BMAI sequence, then for any  irrational $\beta \in (0,1)$  we have 
\begin{align*}
    \lim_{M \to \infty} \E_{m\leq M}a(m) e(m \beta) = 0.
\end{align*}
\begin{proof}
By \cref{orthogonality} we have a universal constant $C>0$ such that for $N,S \in \N $
\begin{equation}\label{Florian-inequality}
    | \E_{n \leq N}  a(n) e(n \beta) | \leq C \sum_{l = 1}^{\frac{\log S}{\log 2}} \frac{1}{l \log \log S} \left(\max_{p,q \in [2^l,2^{l+1}]} \left\vert\E_{n \leq \frac{N}{\max \{p,q\}}} a(pn) e(pn \beta) \overline{a(qn) e( qn \beta)}\right\vert\right)^{1/2} + \epsilon(S,N),
\end{equation}
with $\limsup_{ S \to \infty } \limsup_{N \to \infty } \epsilon(S,N) = 0$. Notice that 
\begin{align*}
        &\lim_{M \to \infty} \left\vert\E_{ m \leq M} a(pm) \overline{a(qm)} e( \beta (p - q ) m ) \right\vert     \\
         \leq &\lim_{M \to \infty}  \E_{ m \leq M} |( a(pm) \overline{a(qm)}  - \theta_p \overline{\theta_q} )e( \beta (p - q ) m ) | + \lim_{M \to \infty} |\theta_p \overline{\theta_q}|\cdot |\E_{ m \leq M} e( \beta (p - q ) m ) |.
\end{align*}
By lemma \ref{Lemma_BMAI}, we have that the first term of the sum is $o_{p\to\infty}(1)+o_{q\to\infty}(1)$. And as $\beta$ is irrational, we  have that $\beta (p-q)  $ is not an integer and so  $\lim_{M \to \infty} |\E_{ m \leq M} e( \beta (p - q ) m )| = 0$. Therefore, we have
\begin{equation}\label{Florian-inequality}
    \limsup_N | \E_{n \leq N}  a(n) e(n \beta) | \leq C \sum_{l = 1}^{\frac{\log S}{\log 2}} \frac{1}{l \log \log S} \left(\max_{p,q \in [2^l,2^{l+1}]} o_{p\to\infty}(1)+o_{q\to\infty}(1)\right)^{1/2} + \limsup_N \epsilon(S,N).
\end{equation}
Take $L$ big enough such that $(\max_{p,q \in [2^l,2^{l+1}]} o_{p\to\infty}(1)+o_{q\to\infty}(1))^{1/2} \leq\epsilon$, $\forall l\geq L$. Then, taking $\limsup_{S}$ in \cref{Florian-inequality}
\begin{align*}
     \limsup_N | \E_{n \leq N}  a(n) e(n \beta) | &\leq C \limsup_{S}\sum_{l =L}^{\frac{\log S}{\log 2}} \frac{1}{l \log \log S} \left(\max_{p,q \in [2^l,2^{l+1}]} o_{p\to\infty}(1)+o_{q\to\infty}(1)\right)^{1/2}\\
     &\leq  C \epsilon \limsup_{S} \frac{1}{ \log \log S}\sum_{l =L}^{\frac{\log S}{\log 2}} \frac{1}{l}= C \epsilon.
\end{align*}
 Taking $\epsilon\to 0$, we conclude that
$$\lim_N \E_{n \leq N}  a(n) e(n \beta)=0. $$
\end{proof}
\end{Lemma}

\begin{Proposition} \label{proposition_baby_case}
Suppose $a : \N \to \C$ is a BMAI sequence such that  $\lim_{N \to \infty} \frac{1}{N} \sum_{ n = 1}^N a(n) = c$ with $c \in \C$. Then for any irrational $\alpha >1$, we have
\begin{align*}
    \lim_{N \to \infty} \frac{1}{N} \sum_{ n = 1}^N a(\lfloor\alpha n \rfloor ) = c.
\end{align*}
\end{Proposition}

\begin{proof}
First notice that 
\begin{align*}
     \frac{1}{N} \sum_{ n = 1}^N a(\lfloor \alpha n \rfloor ) = \alpha \frac{1}{ \alpha N} \sum_{ m = 1}^{ \alpha N}  a(m) \one_{\lfloor\alpha \N\rfloor}(m).
\end{align*}
And so to prove that  $\lim_{N \to \infty} \frac{1}{N} \sum_{ n = 1}^N a(\lfloor \alpha n \rfloor ) = c$, it suffices to show that 
\begin{align}
    \lim_{M \to \infty} \frac{1}{M} \sum_{ m = 1}^{ M }  a(m) \one_{\lfloor\alpha \N\rfloor}(m) = \alpha ^{-1} c.
\end{align}
By \cref{IndIsAP}, $\one_{\lfloor \alpha \N \rfloor}(m)$ is an almost periodic function such that for $\epsilon>0$, there is $(c_i)_{i=0,\cdots,k} \in \C^k$ and irrationals $(\beta_i)_{i=1,\cdots,k} \in (0,1)^k$ such that
\begin{align*}
   \left\Vert  \one_{\lfloor\alpha \N\rfloor}(m) -\left(c_0 + \sum_{i =1}^{k} c_i e(m \beta_i)\right)\right\Vert _1<\epsilon.
\end{align*}
In addition, as $d([\alpha N]) = \alpha^{-1}$ and for each $\beta\in(0,1)$, $\E_{n\leq N} e(m\beta_i)\to0$, we have that $c_0 = \alpha^{-1}$. \\

Hence, we have
\begin{align*}
    &\limsup_M|\E_{m\leq M} a(m)\one_{\lfloor\alpha \N\rfloor}(m)-\alpha^{-1}c |\\
    & \leq\Vert a\Vert _\infty \left\Vert  \one_{[\alpha \N]}(m) -\left(c_0 + \sum_{i =1}^{k} c_i e(m \beta_i)\right)\right\Vert _1\\
    &\tab+ \limsup_M \left\vert\E_{n\leq M} a(m)\left(c_0 + \sum_{i =1}^k c_i e(m \beta_i)\right) -\alpha^{-1}c \right\vert \\
    & \leq \Vert a\Vert _\infty \epsilon + \underbrace{\limsup_M |\E_{n\leq M} a(m)\alpha^{-1} - c\alpha^{-1}|}_{=0} + \sum_{i =1}^k|c_i|\limsup_M |\E_{n\leq M} a(m)e(m \beta_i))  |.
\end{align*}
By the assumptions on $a$, $\lim_{M \to \infty} \frac{1}{M} \sum_{ m = 1}^{ M }  a(m)\alpha^{-1} = \alpha^{-1}c$. By \cref{Lemma_BMAI_and_exp}, the last term of the sum tends to $0$, giving
$$ \limsup_M|\E_{n\leq M} a(m)\one_{\lfloor\alpha \N\rfloor}(m)-\alpha^{-1}c |\leq \Vert a\Vert _\infty \epsilon.$$
We conclude by taking $\epsilon \searrow 0$.
\end{proof}

\begin{Remark}
Note that the condition of the above proposition can be weakened. Indeed, we can remove the conditions of $a$ being BMAI and only ask for a bounded $a : \N \to \C$ which is orthogonal to $e(n\beta)$ for $\beta$ irrational, i.e. $\lim_{M \to \infty} \frac{1}{M} \sum_{ m = 1}^{ M }  a(m) e(m \beta) = 0$. 
\end{Remark}
\begin{Remark}
By von Mangoldt \cite[p. 852]{vonMangoldt97} and Landau \cite[pp. 571–572,
620-621]{Landau09b}, the prime number theorem is equivalent to 
$$\lim_{N\to \infty} \E_{n\leq N} \lambda(n)=0. $$
\end{Remark}

\begin{Corollary}
For any irrational $\alpha > 1$, we have  
\begin{align*}
   \lim_{N \to \infty}  \frac{1}{N} \sum_{n=1}^N \lambda (\lfloor \alpha n\rfloor) = 0. 
\end{align*}
\end{Corollary}
\begin{proof}
Using that
 $$ \lim_{N \to \infty}  \frac{1}{N} \sum_{n=1}^N \lambda (n)= 0,$$
 and the fact that $\lambda$ is BMAI, by \cref{proposition_baby_case}, we obtain the desired result.
\end{proof}

\subsection{Main Theorem}
We now state the main technical result, whose proof will come later. The main theorem follows from this result using \cref{1->2}.

\begin{Theorem}\label{technicalTheorem}
Let $a:\N\to \C$ be a bounded sequence and $b: \N \to  \mathbb{S}^1$ be a BMAI sequence such that: 
\begin{enumerate}
    \item  $\displaystyle \lim_{H \to \infty} \lim_{N \to \infty}  \E_{n\leq N}| \E_{h\leq H}a(n+h)| = 0,$
    \item for all $\beta\in\R\setminus \Q$, $p \neq q$ primes and $ (i,j) \in \{0, \cdots p-1\} \times \{0, \cdots q-1\} \setminus (0,0) $,
    $$\lim_{N \to \infty}  \E_{n\leq N}a(pn + i)\overline{a(qn+j)} e(n \beta) = 0 ~~~ ( \text{resp. }\E^{\log}_{n\leq N} \bullet).$$
\end{enumerate}
Then  for any irrational $\alpha \in \R_+$,

\begin{equation*}
\lim_{N \to \infty} \E_{n\leq N} b(n) a(\lfloor \alpha n \rfloor) = 0 ~~~( \text{resp. } \E^{\log}_{n\leq N}\bullet ).
\end{equation*}
\end{Theorem}
We now  use \cref{technicalTheorem} to prove \cref{Main_Theorem}.

\begin{proof}[Proof of \cref{Main_Theorem} using \cref{technicalTheorem}]
We just need to prove condition 2 from \cref{technicalTheorem}. Let  $\beta\in\R\setminus \Q$, $p \neq q$ primes and $ (i,j) \in \{0, \cdots p-1\} \times \{0, \cdots q-1\} \setminus (0,0)$. If $iq=pj$, then $p| i$, which implies that $i=0$ and, therefore, $j=0$. Nevertheless, this contradicts that $(i,j)\neq (0,0)$ and, in consequence, $iq-pj\neq 0$. Finally, using Lemma \ref{1->2} we immediately obtain condition 2.
\end{proof}
The following example illustrates the necessity of the first condition in both theorems.
\begin{Example}
Consider $a\equiv b\equiv 1$. Clearly both functions are BMAI with values in $\mathbb{S}^1$, and satisfy condition 2 of \cref{Main_Theorem} by
$$\Big|\frac{1}{N}\sum_{n=1}^N e(n\beta)\Big|\leq \frac{1}{N}  \cdot \frac{2}{|1-e(\beta)|}\to 0. $$

However, $a$ clearly doesn't satisfy condition $1$ as the associated limit is $1$.
\end{Example}

In contrast, the BMAI condition can be relaxed. \cref{cor-notbmai} below is an example of a function which is not BMAI but satisfies the conditions of \cref{Main_Theorem}.

\begin{Corollary}\label{cor-notbmai}
For an irrational $\alpha$ and a real number $c$,
\begin{align}\label{lambda-e-limit}
   \lim_{N \to \infty }\frac{1}{N} \sum_{n=1}^N \lambda(n) e(c\lfloor \alpha n \rfloor )  = 0.
\end{align}
Furthermore, for any almost periodic function $f : \N \to \C$,
\begin{align}\label{lambda-f-limit}
   \lim_{N \to \infty }\frac{1}{N} \sum_{n=1}^N \lambda(n) f(\lfloor \alpha n \rfloor )  = 0.
\end{align}

\end{Corollary}
\begin{Remark} Similar results have been obtained in other contexts such as \cite{GT12b}.
\end{Remark}
\begin{proof}
If $c \in \Z$, then \eqref{lambda-e-limit} reduces to \cref{Chowlas-conjecture} (Chowla's conjecture) for the case $k=1$, which is well known to be equivalent to the prime number theorem.\\

Suppose $c \in \R\setminus \Z$, then
\begin{align*}
   \frac{1}{N} \sum_{n=1}^N | \frac{1}{H} \sum_{h =1 }^H e(cn+ch)|  =  \frac{1}{N} \sum_{n=1}^N | \frac{1}{H} \sum_{h =1 }^H e(ch)| = | \frac{1}{H} \sum_{h =1 }^H e(ch)|
\end{align*}
as $c \notin \Z $, $ \lim_{H \to \infty}  | \frac{1}{H} \sum_{h =1 }^H e(ch)| = 0$
and $e(cn)$ satisfies condition 1 and 2 of \cref{Main_Theorem}. Since $\lambda$ is a BMAI sequence, the result follows from \cref{Main_Theorem}. \\

As for \eqref{lambda-f-limit}, note that one can approximate almost periodic functions by a finite sum of $e(cn)$ functions. 
\end{proof}

The following lemma will be useful in proving \cref{technicalTheorem}.

\begin{Lemma}\label{previousLemma}
For any bounded $a: \N \to \C$ that satisfies condition $1$ of \cref{Main_Theorem},
\begin{align*}
        \lim_{H \to \infty} \limsup_{N \to \infty}  \E_{n\leq N} |  \E_{h\leq H}  a(Hn+h)| = 0 ~~~ ( \text{resp. }\E^{\log}_{n\leq N}|\E_{h\leq H}\bullet | ).
\end{align*}

\end{Lemma}
\begin{proof}
Consider the function
\[
\gamma_H(n)= \left\vert  \frac{1}{H} \sum_{h = 1}^H a(n+h) \right\vert,
\]
and let $0<\epsilon<1$. By assumption, we can find an $\tilde{H}\in \N$ such that for every $H\geq \tilde{H}$, $$\limsup_{N}\E_{n\leq N} \gamma_H(n)\leq \epsilon^3. $$
In the following, we will omit the subscript $H$ in $\gamma_H$ and just write $\gamma$. All statements will hold for any $H \geq \tilde{H}$. Let $E=\{n\in \N \mid \gamma(n)\leq \epsilon\}$, and note that
\begin{equation}\label{eq1}
   \E_{n\leq N} \gamma(n) \geq \E_{n\leq N} \gamma(n) \one_{E^c}(n) \geq \E_{n\leq N} \epsilon \one_{E^c}(n).
\end{equation}
Taking $\limsup$ in \cref{eq1},
$$\epsilon^3 \geq \limsup_N\E_{n\leq N} \gamma(n) \geq \epsilon d(E^c) \Longrightarrow \epsilon^2 \geq d(E^c).$$
Let
$$D=\{ n\in \N \mid  \frac{1}{H} |[Hn,H(n+1) )\cap E| \geq (1-\epsilon)\}=\{ n\in \N \mid  \frac{1}{H} |[Hn,H(n+1) )\cap E^c| \leq \epsilon\} $$ be the set of numbers $n$ for which the $H$ points after starting points $Hn$ are covered by $E$ up to $\epsilon$. Note that 
\begin{align*}
\E_{n\leq (N+1)H-1} \one_{E^c} &= \frac{NH}{(N+1)H-1}\cdot \frac{1}{NH} \sum_{n\leq N} |[Hn,H(n+1))\cap E^c| \\
&\geq  \frac{NH}{(N+1)H-1}\cdot \epsilon \frac{1}{N} \sum_{n\leq N} \one_{D^c}(n).
\end{align*}
Taking $\limsup$, we conclude that
$$\epsilon^2\geq \epsilon d(D^c) \Longrightarrow \epsilon\geq d(D^c). $$
In other words, since $E$ has high density, ``good" starting points are also ``dense".
If $n\in D$ then $Hn$ is at most at $\epsilon H$ distance of an element of $m\in E$. Using that $\gamma(m) \leq \epsilon \Rightarrow \gamma(m-k)$ $\leq \epsilon+\frac{2k\| a\|_\infty }{H}$ for $0 \leq k \leq \epsilon H$ we conclude that $\gamma(Hn)\leq \epsilon +2\epsilon\| a\|_\infty$ for every $n\in D$.
Finally, we note that
$$\E_{n\leq N} \gamma(nH)=\E_{n\leq N} \gamma(nH)\one_D(n)+\E_{n\leq N} \gamma(nH)\one_{D^c}(n)\leq  \epsilon +2\epsilon\| a\|_\infty+ \Vert a\Vert _\infty\E_{n\leq N}\one_{D^c}(n) $$
$$\implies \limsup_N \E_{n\leq N} \gamma(nH)\leq \epsilon +2\epsilon\| a\|_\infty+\Vert a\Vert _\infty d(D^c)\leq \epsilon(1+3 \Vert a\Vert _\infty).$$
As this is valid for every $H\geq \tilde{H}$, we conclude that
\begin{equation} \label{equation_cond_1_modified}
\lim_{H\to \infty} \limsup_N \E_{n\leq N} \gamma(nH)=0.
\end{equation}

For the logarithmic version, we use \cref{equation_cond_1_modified} combined with \cref{Lemma_zero_Ces_Log}, taking $S = \N $ and $f(n,T) = \gamma_T(n)$.

\end{proof}
\begin{proof}[Proof of \cref{technicalTheorem}]
For simplicity, we restrict our attention to Cesaró averages, as the proof for logarithmic averages is analogous and all the necessary lemmas used in the proof hold for both types of averages. We will apply  \cref{orthogonality} to $f(n) = b(n) a(\lfloor \alpha n \rfloor)$ to get
\begin{align*}
    &|\E_{ n \leq N} b(n) a(\lfloor \alpha n \rfloor) | \\
    &\leq C \frac{1}{\log \log S} \sum_{l = 1}^{\frac{\log S}{ \log 2}} \frac{1}{l} (\max_{p,q \in [2^l, 2^{l+1}) } |\E_{n \leq \frac{N}{\max\{p,q\}}} b(pn)\overline{b(qn)}a(\lfloor p \alpha n \rfloor)\overline{a(\lfloor q \alpha n \rfloor)}|)^{\frac{1}{2}} +\epsilon(S,N),
\end{align*}
 with $\limsup_{ S \to \infty } \limsup_{N \to \infty } \epsilon(S,N) = 0$.
Similarly as in the proof of \cref{proposition_baby_case}, we take limits and use \cref{Lemma_BMAI} to get
\begin{equation}
\begin{split}
    &\limsup_{N \to \infty} |\E_{ n \leq N} b(n) a(\lfloor \alpha n \rfloor) | \leq  \\
    &\limsup_{S \to \infty} \limsup_{N \to \infty}  \frac{C\Vert b\Vert_{\ell^\infty}}{\log \log S} \sum_{l = 1}^{\frac{\log S}{ \log 2}} \frac{1}{l} (\max_{p,q \in [2^l, 2^{l+1}) } |\theta_p \overline{\theta_q}| |\E_{n \leq \frac{N}{\max \{p,q\}}} |b(n)|^2 a(\lfloor p \alpha n \rfloor)\overline{a(\lfloor q \alpha n \rfloor)}|)^{\frac{1}{2}}.
    \end{split}
\end{equation}

Moreover, as $|\theta_p \overline{\theta_q}| \leq C_\theta$ for all $p \neq q $, we get that 
\begin{equation} \label{Katai_Main_Theorem}
\begin{split}
    &\limsup_{N \to \infty} |\E_{ n \leq N} b(n) a(\lfloor \alpha n \rfloor) | \leq \\
    &\limsup_{S \to \infty} \limsup_{N \to \infty}  \frac{C\Vert b\Vert _{\ell^\infty}C_\theta}{\log \log S} \sum_{l = 1}^{\frac{\log S}{ \log 2}} \frac{1}{l} (\max_{p,q \in [2^l, 2^{l+1}) }  |\E_{n \leq \frac{N}{\max\{p,q\}}} a(\lfloor p \alpha n \rfloor)\overline{a(\lfloor q \alpha n \rfloor)}|)^{\frac{1}{2}}.
    \end{split}
\end{equation}
 It is thus sufficient to show
\begin{equation}\label{Eq3.11}
        \limsup_{S \to \infty} \limsup_{N \to \infty}  \frac{1}{\log \log S} \sum_{l = 1}^{\frac{\log S}{ \log 2}} \frac{1}{l} (\max_{p,q \in [2^l, 2^{l+1}) }  |\E_{n \leq \frac{N}{\max \{p,q\}}} a(\lfloor p \alpha n \rfloor)\overline{a(\lfloor q \alpha n \rfloor)}|)^{\frac{1}{2}} = 0.
\end{equation}
To begin, note that  
\begin{align}\label{bigp}
       &\frac{1}{\log \log S} \sum_{l = 1}^{\log \frac{1}{\alpha}} \frac{1}{l} (\max_{p,q \in [2^l, 2^{l+1}) }  |\E_{n \leq \frac{N}{\max \{p,q\}}} a(\lfloor p \alpha n \rfloor)\overline{a(\lfloor q \alpha n \rfloor)}|)^{\frac{1}{2}}\\ 
       &\leq  \frac{1}{\log \log S} \sum_{l = 1}^{\log \frac{1}{\alpha}} \frac{1}{l} (\max_{p,q \in [2^l, 2^{l+1}) }  |\E_{n \leq \frac{N}{\max \{p,q\}}} \|a\|_\infty|)^{\frac{1}{2}}\\
       &\leq \frac{\|a\|_\infty^2}{\log \log S} \sum_{l = 1}^{\log \frac{1}{\alpha}} \frac{1}{l} \\
       &\leq O_\alpha(\frac{1}{\log \log S}), 
\end{align}
which tends to $0$ as $S$ goes to infinity. 
Thus, it remains to investigate the case where $p,q > \frac{1}{\alpha}$. 
Let $i(p,n) \coloneqq \lfloor p \alpha n \rfloor - p \lfloor \alpha n \rfloor$, thus 
$$a(\lfloor p \alpha n \rfloor)\overline{a(\lfloor q \alpha n \rfloor)} =
a(p \lfloor \alpha n \rfloor + i(p,n))\overline{a(q \lfloor \alpha n \rfloor + i(q,n))}.$$ 
Notice that $i(p,n)$ only takes values in $\{ 0, \cdots, p-1 \}$. Let
$$D_{p,i} \coloneqq \{n\in \lfloor \alpha \N\rfloor: n=\lfloor m\alpha \rfloor \text{ and } \lfloor p \alpha m\rfloor = p\lfloor \alpha m \rfloor +i\}. $$ $D_{p,0}, \dots, D_{p,p-1}$ form a partition of $\lfloor \alpha \N \rfloor$.
We have 
\begin{align}\label{Eq1.9}
    \E_{n \leq \frac{N}{\max \{p,q\}}} a(\lfloor p \alpha n \rfloor )\overline{a(\lfloor q \alpha n \rfloor)} 
    &=\alpha \E_{n \leq \frac{\alpha N}{\max \{p,q\}}} \sum_{i(p)} \sum_{j(q)} a(p n  + i)\overline{a(q n + j)}\one_{D_{p,i}}(n) \one_{D_{q,i}}(n)
\end{align}
where $i(p)$ denotes $0\leq i < p$ and similarly for $j(q)$. We consider two cases separately, that of $(i,j) = (0,0)$ and $(i,j) \neq (0,0)$.\\

Firstly, for the case $(i,j)=(0,0)$, since $\alpha$ is irrational, the density of the set $D_{p,0} \cap D_{q,0}$ is $\frac{1}{pq}$. Hence
$$\limsup_{N \to \infty} \E_{n \leq \frac{N}{\max \{p,q\}}} | a(p n  + i)\overline{a(q n + j)} \one_{D_{p,0} \cap D_{q,0}} (n)| \leq  \frac{ \Vert a\Vert _\infty^2}{pq}.$$ 
 Therefore, the contribution of $(i,j) = (0,0)$ in \cref{Eq3.11} is
\begin{align*}
   \leq \frac{\alpha \Vert a\Vert _\infty}{\log \log S} \sum_{l = \log \frac{1}{\alpha}}^{\frac{\log S}{ \log 2}} \frac{1}{l} \max_{p,q \in [2^l, 2^{l+1}) } \frac{1}{\sqrt{pq}} \leq  \frac{\alpha \Vert a\Vert _\infty}{\log \log S} \sum_{l = 1}^{\frac{\log S}{ \log 2}} \frac{1}{l2^l} \to_{S \to \infty } 0 .
\end{align*}
\\
We now consider the case $(i,j) \neq (0,0)$. We want to estimate
\begin{align}
    \limsup_N |\sum_{\stackrel{i(p),j(q)}{(i,j) \neq (0,0)}} \alpha \E_{m \leq \frac{N \alpha}{2^l}} a(p m + i)\overline{a(q m + j)}\one_{D_{p,i} \cap D_{q,i} }(m)|. \label{equation_proof_main_Theorem_split}
\end{align}
For each $(i,j) \neq (0,0)$, $\one_{D_{p,i} \cap D_{q,i} }(m) $ is an almost periodic function that can be extended to a piece-wise continuous function on all of $\R$ with period $\alpha$, using \cref{DpiAP}. In particular, for $\epsilon>0$, by Lemma \ref{APApproximation} there exist $c_1,\cdots, c_K\neq 0$ and $\beta_1,\cdots,\beta_K\in \R\setminus \Q$ such that
\begin{equation*}
    \Vert  \one_{D_{p,i} \cap D_{q,i} }(m)-\Bigl( \frac{1}{pq} + \sum_{k=1}^K c_k e(m \beta_k) \Bigr)\Vert _1<\frac{\epsilon}{\Vert a\Vert _\infty^2}.
\end{equation*}
Using the triangle inequality,
\begin{align*}
    \cref{equation_proof_main_Theorem_split} &\leq
    \limsup_N \alpha \E_{m \leq \frac{N \alpha}{2^l}} | \sideset{}{'} \sum_{i(p),j(q)}  \frac{1}{pq}  a(p m + i)\overline{a(q m + j)}| \\
    &+  \limsup_N   \alpha \sum_{\stackrel{i(p),j(q)}{(i,j) \neq (0,0)}} \sum_{k=1}^K c_k | \E_{m \leq \frac{N \alpha}{2^l}}a(p m + i)\overline{a(q m + j)}e(m \beta_k)| + \epsilon.
\end{align*}
The second $\limsup_N$ is $0$ from our assumptions on $a$. We are thus left with the first term of the sum. Rewriting,
\begin{align}
    \E_{m \leq \frac{N \alpha}{2^l}} |\sum_{\stackrel{i(p),j(q)}{(i,j) \neq (0,0)}} \frac{1}{pq}  a(p m + i)\overline{a(q m + j)}| \leq \E_{m \leq \frac{N \alpha}{2^l}} | \frac{1}{p} \sum_{i = 0}^{p-1}a(p m + i) |\cdot | \frac{1}{q} \sum_{j=0}^{q-1} \overline{a(q m + j)}  |.
\end{align}
Using \cref{previousLemma} we have that 
\begin{align*}
        \lim_{H \to \infty} \limsup_{N \to \infty}  \E_{n\leq N} |  \E_{h\leq H}  a(Hn+h)| = 0,
\end{align*}
and, taking $l$ large enough (and, therefore, $H=p$ large enough) we derive that
\begin{align*}
    &\limsup_{N} \E_{m\leq \frac{N\alpha}{2^l}} | \frac{1}{p} \sum_{i = 0}^{p-1}a(p m + i) | \leq \frac{\epsilon}{\Vert a\Vert _\infty} \\
    &\implies \limsup_{N}\E_{m \leq \frac{N \alpha}{2^l}} | \frac{1}{p} \sum_{i = 0}^{p-1}a(p m + i) |\cdot |\frac{1}{q} \sum_{j=0}^{q-1} \overline{a(q m + j)}  | \leq \epsilon.
\end{align*}
Thus for $l$ large enough 
$$    \limsup_N |\sum_{\stackrel{i(p),j(q)}{(i,j) \neq (0,0)}} \alpha \E_{m \leq \frac{N \alpha}{2^l}} a(p m + i)\overline{a(q m + j)}\one_{D_{p,i} \cap D_{q,i} }(m)|\leq 3 \epsilon. $$

Finally, together with \cref{bigp}, we have that 
$$\limsup_{S \to \infty} \limsup_{N \to \infty}  \frac{1}{\log \log S} \sum_{l = 1}^{\frac{\log S}{ \log 2}} \frac{1}{l} (\max_{p,q \in [2^l, 2^{l+1}) }  |\E_{n \leq \frac{N}{2^l}} a(\lfloor p \alpha n \rfloor)\overline{a(\lfloor q \alpha n \rfloor)}|)^{\frac{1}{2}}\leq 3\epsilon,$$
whence taking $\epsilon\to 0$ gives the desired result.
\end{proof}
\section{Consequences of the main theorem}\label{sec4}
\subsection{\Erdos{}-Kac generalization}
For the following sections, we denote for every $n\in \N$,
 $$B_n\coloneqq\frac{\omega(n)-\log\log{n}}{(\log\log{n})^{1/2}},~~A_n\coloneqq\frac{\omega(\lfloor\alpha n\rfloor)-\log\log{n}}{(\log\log{n})^{1/2}},$$
and for a function $F:\R\to \C$,\\
\begin{enumerate}
    \item $\Psi_F(n) \coloneqq F(B_n)$,
    \item $I(F) \coloneqq \frac{1}{\sqrt{2\pi}} \int\limits_{-\infty}^{\infty} F(x)e^{-\frac{x^2}{2}} dx$.\\
\end{enumerate}

First, we prove that $a=\Psi_G$ satisfies the conditions of theorem \cref{technicalTheorem} whenever $I(G)=0$.

\begin{Lemma}\label{EKConditions}
Let $G : \R \to \C$ be a bounded continuous function such that $I(G) = 0$. Then $\Psi_G$ satisfies the conditions of \cref{Main_Theorem}.
\end{Lemma}
\begin{proof}
First, we show that $\Psi_G$ satisfies the first condition. Proving the following will be sufficient:
\begin{align*}
\lim_{H \to \infty} \lim_{N \to \infty}  \frac{1}{N} \sum_{n = 1}^N |\frac{1}{H} \sum_{h = 1}^H \Psi_G(n+h)|^2 = 0.
\end{align*}
Indeed, by expanding the squares and switching the order of summation, we get 
\begin{align*}
\lim_{H \to \infty} \lim_{N \to \infty}  \frac{1}{N} \sum_{n = 1}^N |\frac{1}{H} \sum_{h = 1}^H \Psi_G(n+h)|^2 = \lim_{H \to \infty} \lim_{N \to \infty} \frac{1}{H^2} \sum_{h = 1}^H \sum_{h' = 1}^H \frac{1}{N} \sum_{n = 1}^N \Psi_G(n+h) \overline{\Psi_G(n+h')}.
\end{align*}
Splitting the sum into the diagonal terms and the non diagonal terms, 
\begin{align*}
\lim_{H \to \infty} \frac{1}{H^2} \sum_{h \neq h'} \lim_{N \to \infty}  \frac{1}{N} \sum_{n = 1}^N \Psi_G(n+h) \overline{\Psi_G(n+h')} + \lim_{H \to \infty} \lim_{N \to \infty}  \frac{1}{H^2} \sum_{h=1}^H\frac{1}{N} \sum_{n = 1}^N |\Psi_G(n+h)|^2.
\end{align*}
 The first sum tends to $0$. Indeed let $\epsilon>0$ and $(E_i)_{i=1}^I$ be Jordan measurable sets, where $E_i \subseteq \R$ for each $i \in I$, and $(c_i)_{i=1}^I \subset \C^I$ such that 
 
 $$||\sum_{i= 1}^Ic_i\one_{E_i}- G||_{L^1(\R)}\leq \sqrt{2\pi \epsilon} ,~~~~~~ ||\sum_{i= 1}^Ic_i\one_{E_i}||_{L^\infty(\R)}\leq ||G||_{L^\infty(\R)}.$$
We may assume $||\one_{[-M,M]} (\sum_{i= 1}^Ic_i\one_{E_i}- G)||_{L^\infty(\R)}\leq \epsilon $, where $M>0$ is such that 
$$|\int_{[-M,M]^c}  e^{-x^2/2} dx|\leq \frac{1}{2||G||_{L^\infty(\R)}}\epsilon. $$ 

 Notice that 
\begin{align*}
     & \limsup_{n\leq N}\E_{n\leq N}  |\Psi_G(n+h) - \sum_{i= 1}^Ic_i\one_{E_i}(B_{n+h})|\\
& \leq\limsup_{N \to \infty}\E_{n\leq N} |\Bigg(\one_{[-M,M]}\Big(G - \sum_{i= 1}^Ic_i\one_{E_i}\Big)\Bigg)(B_{n+h})|\\
&+\limsup_{N \to \infty}\E_{n\leq N} |\Bigg(\one_{[-M,M]^c}\Big(G - \sum_{i= 1}^Ic_i\one_{E_i}\Big)\Bigg)(B_{n+h}) |\\
& \leq ||\one_{[-M,M]} (\sum_{i= 1}^Ic_i\one_{E_i}- G)||_{L^\infty(\R)}+ 2||G||_{L^\infty(\R)} \limsup_{N \to \infty}\E_{n\leq N}\one_{[-M,M]^c}(B_{n+h})\\
&\leq \epsilon + 2||G||_{L^\infty(\R)}|\int_{[-M,M]^c}  e^{-x^2/2} dx|\leq 2\epsilon,
\end{align*}
where we used \cref{tanaka}. Hence 
 \begin{align*}
      \limsup_{N \to \infty}|\E_{n\leq N}\Big(  \Psi_G(n+h) \overline{\Psi_G(n+h')} - \sum_{i,j\in[I]^2}c_i\overline{c_j}\one_{E_i}(B_{n+h}) \one_{E_j}(B_{n+h'})\Big)|\leq 4\epsilon ||G||_\infty  \end{align*}
Using \cref{tanaka} again, we have

 \begin{align*}
   \lim_{N\to \infty}\E_{n\leq N}\sum_{i,j\in[I]^2}c_i\overline{c_j}\one_{E_i}(B_{n+h}) \one_{E_j}(B_{n+h'})&= \frac{1}{2\pi}\sum_{i,j\in[I]^2}c_i\overline{c_j} \int_{E_i\times E_j} e^{-\frac{1}{2} (x_1^2+x_2^2)} dx_1 dx_2\\
   &=  \frac{1}{2\pi}  \sum_{i,j\in[I]^2}c_i\overline{c_j} (\int_{E_i} e^{-x^2/2} dx)(\int_{E_j} e^{-x^2/2} dx) \\
   &= \frac{1}{2\pi}| \int_{-\infty}^\infty   \sum_{i= 1}^Ic_i\one_{E_i} e^{-x^2/2} dx|^2. 
 \end{align*}
 Using this, we derive 
\begin{align*}
    \limsup_{n\leq N}|\E_{n\leq N}  \Psi_G(n+h) \overline{\Psi_G(n+h')}|& \leq 4\epsilon ||G||_\infty +  \lim_{n\to N}|\E_{n\leq N} \sum_{i,j\in[I]^2}c_i\overline{c_j}\one_{E_i}(B_{n+h}) \one_{E_j}(B_{n+h'})  |  \\
    & = 4\epsilon ||G||_\infty +\frac{1}{2\pi} |\int_{-\infty}^\infty   \sum_{i= 1}^Ic_i\one_{E_i}(x) e^{-x^2/2} dx|^2 \\
&\leq 4\epsilon ||G||_\infty  +\frac{1}{2\pi}|  \int_{-\infty}^\infty   (\sum_{i= 1}^Ic_i\one_{E_i}(x)- G(x)) e^{-x^2/2} dx |^2 \\
&\text{ (since $I(G)= 0$) }\\
& \leq  4\epsilon ||G||_\infty  +\frac{1}{2\pi} ||\sum_{i= 1}^Ic_i\one_{E_i}- G||_1^2\\
&\leq  4\epsilon ||G||_\infty + \epsilon.
\end{align*}
Thus as $\epsilon>0$ was arbitrary, we get that $\lim_{N \to \infty}|\E_{n\leq N}  \Psi_G(n+h) \overline{\Psi_G(n+h')}|=0$.\\


For the second sum,  notice that $\frac{1}{N} \sum_{n = 1}^N |\Psi_G(n+h)|^2 \leq C$ for some absolute constant as $G$ is bounded, and thus
$$\lim_{H \to \infty} \lim_{N \to \infty}  \frac{1}{H^2} \sum_{h=1}^H\frac{1}{N} \sum_{n = 1}^N |\Psi_G(n+h)|^2 \leq \lim_{H \to \infty} \lim_{N \to \infty}  \frac{1}{H^2} \sum_{h=1}^H C = \lim_{H \to \infty} \lim_{N \to \infty} \frac{C}{H} = 0.$$
Hence the first condition is satisfied.\\

To verify the second condition, in light of \cref{corput}, it is enough to show that for all $h \in \N$
\begin{align*}
    \lim_{N \to \infty}  \E_{n\leq N} \Psi_G(p(n+h) + i)\overline{\Psi_G(q(n+h)+j)} e((n+h) \beta)\overline{\Psi_G(pn + i)\overline{\Psi_G(qn+j)} e(n \beta)} = 0,
\end{align*}
which is equivalent to showing
\begin{align*}
        \lim_{N \to \infty}  \E_{n\leq N} \Psi_G(pn  + ph +i)\overline{\Psi_G}(qn + qh+j) \overline{\Psi_G}(pn + i)\Psi_G(qn+j) = 0.
\end{align*}
As $p \neq q$, $ (i,j) \in \{0, \cdots p-1\} \times \{0, \cdots q-1\} \setminus (0,0) $ and $h > 0$, we have that $pn +ph +i$, $qn + qh+j$, $pn+i$ and $qn+j$ are polynomials over $\Z$ that are pairwise relatively prime. We can thus apply \cref{tanaka} to get the desired result.
\end{proof}

\begin{Lemma}\label{lemma3.12}
Let $\alpha\in \R$ be a positive irrational real number. Then for every continuous $F: \R \rightarrow \mathbbm{S}^1$ and bounded $G: \R \rightarrow \mathbbm{C}$ with $I(G) = 0$,  
\begin{equation}
    \lim_{N\to \infty}\E_{n\leq N}\Psi_F(n) \Psi_G(\lfloor \alpha n \rfloor) =0.
\end{equation}
\end{Lemma}
\begin{proof}
Let $a(n)=\Psi_G(n)$  and $b(n)=\Psi_F(n)$. By \cref{EKConditions}, $a$ satisfies the conditions of \cref{Main_Theorem} and $b$ is bounded. Hence
$$   \lim_{N\to \infty}\E_{n\leq N}\Psi_F(n) \Psi_G(\lfloor \alpha n \rfloor) =\lim_{N\to \infty} \E_{n\leq N}b(n)a(\lfloor\alpha n\rfloor)=0. $$
\end{proof}

\begin{Lemma}\label{lemma3.13}
Suppose that for every continuous function with compact support $F : \R \to \mathbb{S}^1$, and bounded function $G:\R\to \C$ with $I(G)=0$,
$$\lim_{N\to \infty} \E_{n\leq N}  F(A_n)G(B_n)= 0. $$ 
Then, for all $a,b,c,d\in \R$,  
 $$\lim_{N\to \infty} \E_{n\leq N} \one_{[a,b]}(A_n) \one_{[c,d]}(B_n)=I(\one_{[a,b]}) I(\one_{[c,d]}).$$
\end{Lemma}
\begin{proof}
First, notice that the family of functions
$$\mathcal{A} = \text{span}\Bigl(F \in \mathcal{C}(\R,\mathbb{S}^1) \mid F:\R \to \mathbb{S}^1 \text{ has compact support }\Bigr), $$
is such that 
\begin{itemize}
    \item $\mathcal{A}$ separates points, in other words, for any $x \neq y \in X$, there exists $f \in \mathcal{A}$ such that $f(x) \neq f(y)$,
    \item for every $x \in X$, there exists $f \in \mathcal{A}$ such that $f(x) = 0$,
    \item $\mathcal{A}$ is self-adjoint i.e. if $f \in \mathcal{A}$ then $\overline{f} \in \mathcal{A}$.
\end{itemize}
Then by the Stone-Weierstrass theorem, $\mathcal{A}$ is dense in $C_\infty(\R,\C)$. If $H\in C_\infty(\R,\C)$ and $\epsilon>0$, let $F_1,\cdots,F_r\in \mathcal{A}$ and $\{c_i\}_{i=1}^r\subseteq \C$ such that 
$$\Vert H- \sum_{i=1}^r c_i F_i\Vert _\infty\leq \epsilon. $$
Then, denoting $F_\epsilon=\sum_{i=1}^r c_i F_i$ and $M=\max_{i=1,\cdots,r}|c_i|$,
\begin{align*}
    {|\E_{n\leq N} H(A_n)G(B_n) |}&\leq |\E_{n\leq N} \sum_{i=1}^r c_i F_i(A_n)G(B_n) |+ \E_{n\leq N} | F_\epsilon(A_n)G(B_n)- H(A_n)G(B_n)|\\
    &\leq  M\sum_{i=1}^r |\E_{n\leq N} F_i(A_n)G(B_n) |+ \Vert G\Vert _\infty  \Vert  F_\epsilon- H\Vert _\infty\\
    &\leq M\sum_{i=1}^r |\E_{n\leq N} F_i(A_n)G(B_n) |+ \Vert G\Vert _\infty\epsilon.
\end{align*}
Thus, using the hypothesis,
$$\limsup_{N\to \infty}  |\E_{n\leq N} H(A_n)G(B_n) | \leq \Vert G\Vert _\infty\epsilon, $$
and as this is valid for all $\epsilon>0$, we conclude that 
$$\lim_{N\to \infty} \E_{n\leq N} H(A_n)G(B_n) =0,~~ \forall H\in C_\infty(\R,\C). $$

Now, given $a,b,c,d\in \R$ such that $a<b, c<d$ (otherwise it is trivial) and $\epsilon>0$ such that $\epsilon<(b-a)/2$, we choose $H\in C_\infty(\R)$ such that 
$$H(x)=\begin{cases} 1  &  x\in [a+\epsilon,b-\epsilon],\\
0& x\in (-\infty,a-\epsilon]\cup [b+\epsilon,\infty),
\end{cases} $$
and such that $H$ is monotone in $(a-\epsilon,a+\epsilon)$ and in $(b-\epsilon,b+\epsilon)$. Therefore, we have that
\begin{align*}
 |\E_{n\leq N}(\one_{[a,b]}(A_n )G(B_n))| &\leq |\E_{n\leq N}(\one_{[a,b]}(A_n )G(B_n))-H(A_n)G(B_n)| +|\E_{n\leq N}(H(A_n )G(B_n))|   \\
 & \leq  \Vert G\Vert _\infty\E_{n\leq N}|\one_{[a,b]}(A_n )-H(A_n)| + |\E_{n\leq N}(H(A_n )G(B_n))|\\
 &\leq \Vert G\Vert _\infty\E_{n\leq N}\one_{(a-\epsilon,a+\epsilon)\cup(b-\epsilon,b+\epsilon)}(A_n ) + |\E_{n\leq N}(H(A_n )G(B_n))|.
\end{align*}

Taking $\limsup$ and using \cref{ErdosKacBS} yields 
$$ \limsup_{N\to \infty}|\E_{n\leq N}(\one_{[a,b]}(A_n )G(B_n))| \leq  \Vert G\Vert _\infty I(\one_{(a-\epsilon,a+\epsilon)\cup(b-\epsilon,b+\epsilon)})\leq C\Vert G\Vert _\infty \epsilon,$$
with $C=4/\sqrt{2\pi}$. Therefore, as this is valid for every $\epsilon>0$, we conclude that 
$$\lim_{N\to \infty}\E_{n\leq N}(\one_{[a,b]}(A_n )G(B_n))=0. $$
Taking $G(x)=\one_{[c,d]}(x)-I(\one_{[c,d]})$, we have that 
$$\lim_{N\to \infty}\E_{n\leq N}\one_{[a,b]}(A_n )(\one_{[c,d]}(B_n)-I(\one_{[c,d]}))=0, $$
and thanks to \Erdos{}-Kac (\cref{ErdosKac}) and \cref{ErdosKacBS}, this is equivalent to 
 $$\lim_{N\to \infty} \E_{n\leq N} \one_{[a,b]}(A_n) \one_{[c,d]}(B_n)=I(\one_{[a,b]}) I(\one_{[c,d]}),$$
 concluding the proof.
\end{proof}

\begin{Lemma}\label{technical lemma}
There exists a set $E\subseteq \N$ of density $1$ such that 
$$|A_n-B_{\lfloor \alpha n\rfloor}|\xrightarrow[\stackrel{n\to \infty}{n\in E}]{}0. $$
\end{Lemma}
\begin{proof}
Notice that $ \log\log \lfloor \alpha n\rfloor=    \log\log n+ b_n$, where $(b_n)_n$ is a bounded sequence (in which we used that $\log:(1,\infty)\to \R$ is $1$-Lipschitz given that its derivative is bounded by $1$), and let $C_n\coloneqq \frac{w(\lfloor \alpha n\rfloor)-\log\log n}{(\log\log \lfloor \alpha n\rfloor )^{1/2}}$. Note that
$$|B_{\lfloor \alpha n\rfloor}-C_n|\leq \frac{|b_n|}{\log\log \lfloor \alpha n\rfloor}\xrightarrow[n\to \infty]{}0.  $$
Therefore, it is enough to see that $|A_n-C_n|\to_n 0$ with $n$ taking values in a set of density $1$. Using similar arguments, we have that $(\log\log \lfloor \alpha n\rfloor)^{1/2}=(\log\log n)^{1/2}+b'_n$, where $(b_n')_n$ is a bounded sequence. Let $M=\max(\Vert (b_n)_n\Vert _{\ell^\infty},\Vert (b_n')_n\Vert _{\ell^\infty})$. Hence
\begin{align*}
    |A_n-C_n| &=|b'_n \frac{w(\lfloor \alpha n\rfloor)- \log\log n}{\log\log n + b_n (\log\log n)^{1/2}}| \\
    & \leq M |\frac{1}{1+ b_n/(\log\log n)} |\cdot | \frac{w(\lfloor \alpha n\rfloor)- \log\log n}{\log\log n}|.
\end{align*}
The last term tends to $0$ with $n$ taking values in a set $E$ of density $1$ given by the Hardy–Ramanujan theorem (\cref{Hardy-Ramanujan}). Meanwhile, the first term tends to $M\cdot 1$ as $n\to \infty$, therefore we conclude that 
$$|A_n-C_n|\xrightarrow[\stackrel{n\to \infty}{n\in E}]{}0. $$
By the triangle inequality, we get the desired result
$$|A_n-B_{\lfloor \alpha n\rfloor}|\xrightarrow[\stackrel{n\to \infty}{n\in E}]{}0. $$
\end{proof}
We now restate and prove \cref{ErdosKacGeneral}.

\begin{reptheorem}{ErdosKacGeneral}
Let $\alpha\in \R$ be a positive irrational real number and $a,b,c,d\in \R$. Then
\begin{align*}
\lim_{N\to\infty} \frac{1}{N} \#\{n\leq N : a\leq \frac{\omega(\lfloor\alpha n\rfloor)-\log\log{n}}{(\log\log{n})^{1/2}}\leq b, c\leq \frac{\omega(n)-\log\log{n}}{(\log\log{n})^{1/2}}\leq d\}\\
=\Bigl(\frac{1}{\sqrt{2\pi}}\int_{a}^be^{-t^2/2 }dt\Bigr)\Bigl(\frac{1}{\sqrt{2\pi}}\int_{c}^d e^{-t^2/2 }dt\Bigr).   
\end{align*}

\end{reptheorem}
\begin{proof}
We want to show that 
$$\lim_{N\to \infty} \E_{n\leq N} \one_{[a,b]}(A_n) \one_{[c,d]}(B_n)=I(\one_{[a,b]}) I(\one_{[c,d]}). $$
Using Lemma \ref{lemma3.13}, it is enough to prove that for every continuous and compactly supported $F:\R\to \mathbb{S}^1$ and bounded $G:\R\to \C$ with $I(G)=0$,
\begin{equation}\label{eq3.9} 
   \lim_{N\to \infty} \E_{n\leq N} F(A_n) G(B_n)=0. 
\end{equation}
In order to prove \cref{eq3.9}, we show that
$$ \lim_{N\to \infty} |\E_{n\leq N} F(A_n) G(B_n)-F(B_{\lfloor \alpha n\rfloor})G(B_n) |=0.  $$
Indeed, since $G$ is bounded,
\begin{align}
    |\E_{n\leq N} F(A_n) G(B_n)-F(B_{\lfloor \alpha n\rfloor})G(B_n) |&\leq  \| G \|_\infty \E_{n\leq N} |F(A_n)-F(B_{\lfloor \alpha n\rfloor})|,
\end{align}
and as $F$ is uniformly continuous and  has compact support, it is sufficient to show 
$$|A_n-B_{\lfloor \alpha n\rfloor}|\xrightarrow[n\to \infty]{}0, $$
with $n$ taking values in a set of density $1$, which is given by Lemma \ref{technical lemma}. In this way, it's enough to prove that for every continuous with compact support $F:\R\to \mathbb{S}^1$ and bounded $G:\R\to \C$ with $I(G)=0$ we have that
\begin{equation}\label{eq3.11}
   \lim_{N\to \infty} \E_{n\leq N}\Psi_F(n) \Psi_G(\lfloor \alpha n\rfloor)=0, 
\end{equation}
which follows from \cref{lemma3.12}, concluding the proof.
\end{proof}

\subsection{Logarithmic averages of the Liouville function}

We now cite a relatively recent theorem of Matomaki and Radziwill.
\begin{Theorem}[see {\cite{MR16b}}]\label{Matomaki-Radziwill}
$$\lim_{H\to \infty}\limsup_{N\to \infty} \frac{1}{N}\sum_{n=1}^N\Bigl|\frac{1}{H}\sum_{h=1}^H \lambda(n+h) \Bigr|=0. $$
\end{Theorem}
Note that this is the first condition of our main theorem, \cref{Main_Theorem}, satisfied by the Liouville function. The following lemma from \cite{FH18a} shows that the Liouville function and the \Mobius{} function also satisfy the second condition of \cref{Main_Theorem}.

\begin{Lemma}[see {\cite[Corollary 1.5]{FH18a}}]\label{LiouvilleConditions} 
Let $\beta\in \R$ be irrational, then
$$\lim_{N\to \infty} \E_{n\leq \N}^{\log{}} e(n\beta) \prod_{j=1}^j\lambda(n+h_j)=0 $$
for all $l\in\N$ and $h_1,\cdots,h_l\in \Z$. Moreover, a similar statement holds with the \Mobius{} function $\mu$ in place of $\lambda$.
\end{Lemma}

We now restate and prove \cref{LiouvilleTheorem}.
\begin{reptheorem}{LiouvilleTheorem}
Let $\alpha\in \R$ an irrational positive real number, then 
$$\lim_{N\rightarrow \infty} \E_{n\leq N}^{\log} \lambda(n)\lambda(\lfloor\alpha n\rfloor)=0.$$
\end{reptheorem}
\begin{proof}
Immediate by taking $b=\lambda$ and $a=\lambda$ in Theorem \ref{Main_Theorem} and using \cref{LiouvilleConditions}.
\end{proof}

\section{Open Questions}
While we proved the relative independence of $n$ and $\lfloor\alpha n\rfloor$ under specific conditions, we have yet to show independence of different Beatty sequences. Namely, one may consider $\lfloor\alpha n\rfloor$ and $[\beta n\rfloor$ instead, which leads to the following conjecture. 
\begin{Conjecture}
Let $a: \N \to \C$ and $b: \N \to \mathbb{S}^1$ be a BMAI sequences such that we also have that  \\
\begin{enumerate}
    \item $\lim_{H \to \infty} \lim_{N \to \infty}  \frac{1}{N} \sum_{n = 1}^N | \frac{1}{H} \sum_{h = 1}^H a(n+h)| = 0$ 
    \item if for all $\beta\in\R\setminus \Q$, $h\in \N$ we have that : \\
    $\lim_{N \to \infty}  \E_{n\leq N} a(n )\overline{a(n+h)} e(n \beta) = 0  ~~~$ ( resp $\E^{\log}_{n\leq N}$ ) 
\end{enumerate}
Then  for any irrational $\alpha, \beta \in \R$ such that $\alpha$ and $\beta$ are rationally independent, 
\begin{equation*}
\lim_{N \to \infty} \E_{n\leq N} b(\lfloor \beta n\rfloor) a(\lfloor \alpha n \rfloor) = 0. ~~~( \text{resp } \E^{\log}_{n\leq N} ) 
\end{equation*}
\end{Conjecture}
Once we have pairwise independence of Beatty sequences, we can then consider a higher-order version. For example, the following conjecture offers a higher-order generalization of \cref{LiouvilleTheorem}. 
\begin{Conjecture}\label{Conjecture-2}
For $\{\alpha_i\}_{i=2}^\ell\in \R\setminus \Q$ rationally independent, we have
\begin{equation*}
\lim_{N \to \infty} \E_{n\leq N}^{\log{}} \lambda(n)\lambda(\lfloor\alpha_2 n\rfloor)...\lambda(\lfloor\alpha_\ell n\rfloor) = 0.
\end{equation*}

A more recent work of Joni Ter\"{a}v\"{a}inen and Aled Walker (\cite{teravainen2023bohr}) proves a generalization of \cref{LiouvilleTheorem}, namely for $\alpha_1,\alpha_2>0$ and $\beta_1,\beta_2\in \R$ such that $\alpha_1/\alpha_2$ is irrational,
$$\lim_{N\to \infty} \E_{n\leq N}^{log} \lambda([\alpha_1 n+\beta_1])\lambda([\alpha_2 n+\beta_2])=0. $$


\end{Conjecture}

\small{
\bibliographystyle{aomalphanomr}
\bibliography{bibTeXlibrary}
\hfill
\textsc{
\hspace{0.6cm}(David Crnčević) École Polytechnique, Institut Polytechnique de Paris, Rte de Saclay, Palaiseau 91128, France}

    \textit{E-mail address:} \texttt{david.crncevic@polytechnique.edu}
\hfill \break 

\textsc{
(Felipe Hernández) Department of Mathematics, École polytechnique fédérale de Lausanne, Rte Cantonale, 1015 Lausanne, Switzerland}

    \textit{E-mail address:} \texttt{felipe.hernandezcastro@epfl.ch}
\hfill \break 

\textsc{
(Kevin Rizk) Department of Mathematics, Stanford University, 450 Jane Stanford Way, Stanford, CA 94305}

    \textit{E-mail address:} \texttt{kevin.b.rizk@gmail.com}
\hfill \break 

\textsc{
(Khunpob Sereesuchart) Department of Mathematics, University of California, Los Angeles, 520 Portola Plaza, Los Angeles, CA 90095}

    \textit{E-mail address:} \texttt{ksereesu2417@math.ucla.edu}
\hfill \break 

\textsc{
(Ran Tao) Department of Mathematical Sciences, Carnegie Mellon University, 5000 Forbes Ave, Pittsburgh, PA 15213}

    \textit{E-mail address:} \texttt{rant2@andrew.cmu.edu}
}

\end{document}